\documentclass[12pt, reqno]{amsart}
\usepackage{amsmath, amsthm, amssymb, stmaryrd}
\usepackage{hyperref}
\usepackage{amsfonts}
\usepackage{enumerate}
\usepackage{url}
\usepackage{comment}
\usepackage{xcolor}

\newtheorem{thm}{Theorem}[section]
\newtheorem{lemma}[thm]{Lemma}
\newtheorem{prop}[thm]{Proposition}
\newtheorem{cor}[thm]{Corollary}

\newtheorem*{claim*}{Claim}

\theoremstyle{definition}
\newtheorem{defn}[thm]{Definition}

\theoremstyle{remark}
\newtheorem{remark}[thm]{Remark}

\numberwithin{equation}{section}

\newcommand{\I}{[0,1)}

\newcommand{\N}{\mathbb{N}}
\newcommand{\what}[1]{\widehat{#1}}

\newcommand{\hidden}[1]{}

\newcommand{\ca}[1]{\mathcal{#1}}

\newcommand{\mmod}[1]{\,\,\mathrm{mod}\,\,#1}

\def\alp{{\alpha}} 
\def\bet{{\beta}}  
\def\gam{{\gamma}} 
\def\del{{\delta}} 

\def\tet{{\theta}}  
\def\kap{{\kappa}}
\def\lam{{\lambda}} \def\Lam{{\Lambda}}

\def\sig{{\sigma}}

\def\ome{{\omega}}  
\def\eps{\varepsilon}
\def\implies{\Rightarrow}

\def\le{\leqslant} \def\ge{\geqslant}
\def\ls{\leqslant} \def\gs{\geqslant}
\def\leq{\leqslant}
\def\geq{\geqslant}

\def\d{{\,{\rm d}}}

\def \sig{{\sigma}}

\def \bE {\mathbb E}

\def \bN {\mathbb N}
\def \bP {\mathbb P}
\def \bQ {\mathbb Q}
\def \bR {\mathbb R}
\def \bZ {\mathbb Z}

\def \ba {\mathbf a}
\def \bb {\mathbf b}

\def \bk {\mathbf k}

\def \bx {\mathbf x}

\def \bxi {{\boldsymbol{\xi}}}

\def \fA {\mathfrak A}

\def \fJ {\mathfrak J}

\def \cA {\mathcal A}
\def \cB {\mathcal B}
\def \cC {\mathcal C}
\def \cD {\mathcal D}
\def \cE {\mathcal E}
\def \cF {\mathcal F}
\def \cG {\mathcal G}

\def \cI {\mathcal I}
\def \cJ {\mathcal J}
\def \cK {\mathcal K}

\def \cM {\mathcal M}
\def \cN {\mathcal N}
\def \cO {\mathcal O}

\def \cR {\mathcal R}
\def \cS {\mathcal S}

\def \cV {\mathcal V}
\def \cW {\mathcal W}

\def \dim {\mathrm{dim}}

\def \sinc {\mathrm{sinc}}

\def \Bad {{\mathrm{\mathbf{Bad}}}}

\def \dimH {\dim_{\mathrm{H}}}
\def \diam {{\mathrm{diam}}}
\def \dd {\mathrm d}

\usepackage{tikz}
\usetikzlibrary{shapes.multipart}
\usetikzlibrary{shapes,arrows}

\begin{document}
	\title[Inhomogeneous Kaufman and Littlewood]{Inhomogeneous Kaufman measures and diophantine approximation}
	\subjclass[2010]{}
	\author{Sam Chow \and Agamemnon Zafeiropoulos \and Evgeniy Zorin}
	
	\address{Mathematics Institute, Zeeman Building, University of Warwick, Coventry CV4 7AL, United Kingdom}
	\email{sam.chow@warwick.ac.uk}
	
	\address{Department of Mathematics, Technion, Haifa, Israel}
	\email{agamemn@campus.technion.ac.il}
	
	\address{Department of Mathematics, University of York,
		Heslington, York YO10 5DD, UK}
	\email{evgeniy.zorin@york.ac.uk}
	
	\begin{abstract}
		We introduce an inhomogeneous variant of Kaufman's measure, with applications to diophantine approximation.
		In particular, we make progress towards a problem related to Littlewood's conjecture.
	\end{abstract}
	
	\maketitle
	
	\setcounter{tocdepth}{1}
	\tableofcontents
	
	\section{Introduction}
	
	\subsection{Kaufman's measures}
	
	It follows from Dirichlet's approximation theorem that if $\alp \in \bR$ then there exist infinitely many $n \in \bN$ such that
	\[
	n \| n \alp \| \le 1.
	\]
	The set of \emph{badly approximable} numbers comprises reals for which the right-hand side cannot be replaced by $o(1)$, i.e.
	\[
	\Bad = \{ \alp \in \bR: 
	\inf_{n \in \bN}
	n \| n \alp \| > 0 \}.
	\]
	For any $s \in (0,1)$, Kaufman \cite{Kau1980} constructed a probability measure $\mu$ supported on $\Bad$ with the following two properties that are important for applications:
	\begin{enumerate}[(i)]
		\item (Large Frostman dimension)
		For any interval $I \subseteq [0,1]$, we have
		\[
		\mu(I) \ll_s \lam(I)^s,
		\]
		where $\lam$ denotes Lebesgue measure.
		\item (Polynomial Fourier decay)
		We have
		\[
		\widehat \mu(\xi) \ll 
		(1 + |\xi|)^{-7/10^4}
		\qquad (\xi \in \bR).
		\]
	\end{enumerate}
	Kaufman's construction was generalised by Queff\'elec and Ramar\'e \cite{QR2003}, and further by Jordan and Sahlsten \cite{JS2016}. Measures with Fourier decay are called 
	\emph{Rajchman measures}, and are independently interesting.
	
	In \cite{PVZZ2022}, it is conjectured that if $\gam \in \bR$ then the set
	\[
	\Bad_\gam = \{ \alp \in \bR:
	\inf_{n \in \bN} n
	\| n \alp - \gam \| > 0 \}
	\]
	supports a probability measure $\mu$ with the same two essential properties. Our main result makes progress towards this conjecture. Define
	\[
	\cB = \{ (\alp,\gam) \in \bR^2:
	\alp \in \Bad \cap \Bad_\gam \}.
	\]
	
	\begin{thm} \label{MainThm}
		There exists $\eta > 0$ such that if
		$s \in (0,2)$ then there exists a probability measure $\nu$ supported on $\cB$ with the following two properties:
		\begin{enumerate}[(i)]
			\item \label{MainThm_point_Frostman_dimension} (Large Frostman dimension) 
			For any ball $I$ of radius 
			$r(I)$, we have
			\begin{equation}
				\label{Frostman}
				\nu(I) \ll_s r(I)^s.
			\end{equation}
			\item \label{MainThm_point_decay} (Polynomial Fourier decay)
			For any $\bxi \in \bZ^2$, 
			we have
			\[
			\widehat \nu(\bxi) \ll_\eta
			(1 + \| \bxi \|_\infty)^{-\eta}.
			\]
		\end{enumerate}
	\end{thm}
	\begin{remark} \label{rem_concrete_eta}
		It will be seen in the proof of Theorem~\ref{MainThm} that we can take $\eta$ to be any positive constant smaller than $1/30$. See Remark~\ref{rem_justification_concrete_eta} following the proof of Theorem~\ref{MainThm}.
	\end{remark}
	
	\subsection{Lacunary approximation}
	
	We presently state a special case of
	\cite[Corollary 1]{PVZZ2022}.
	
	\begin{thm} [Pollington--Velani--Zafeiropoulos--Zorin, 2022]
		Fix a real number $A > 2$, and let $\mu$ be a probability measure supported on $[0,1]$ such that
		\[
		\widehat \mu(\xi) \ll_A 
		(\log |\xi|)^{-A}
		\qquad (|\xi| \to \infty).
		\]
		Let $n_1, n_2, \ldots$ be a lacunary sequence of positive integers, let $\del \in \bR$, let 
		$\psi: \bN \to [0,1]$, and write
		$\cW_\del(\psi)$ for the set of 
		$\bet \in [0,1]$ such that
		\[
		\| n_t \bet - \del \| < \psi(t)
		\]
		has infinitely many solutions $t \in \bN$. Let $\eps > 0$. Then
		\[
		\mu(\cW_\del(\psi)) =
		\begin{cases}
			0, &\text{if } \displaystyle
			\sum_{t=1}^\infty
			\psi(t) < \infty \\ \\
			1, &\text{if } \displaystyle
			\sum_{t=1}^\infty
			\psi(t) = \infty.
		\end{cases}
		\]
	\end{thm}
	
	We establish the following two-dimensional analogue.
	
	\begin{thm} \label{LacApprox}
		Let $n_1, n_2, \ldots$ be a lacunary sequence of positive integers,
		let $A > 2$ be sufficiently large, and let $\nu$ be a probability measure supported on $[0,1]^2$ such that
		\begin{equation}
			\label{LogDecay}
			\widehat \nu(\bxi) \ll_A 
			(\log \| \bxi \|_\infty)^{-A}
			\qquad 
			(\|\bxi \|_\infty \to \infty).
		\end{equation}
		Let 
		$\psi: \bN \to [0,1]$, and write
		$\cW(\psi)$ for the set of 
		$(\bet,\del) \in [0,1]^2$ such that
		\begin{equation}
			\label{LacIneq}
			\| n_t \bet - \del \| \le \psi(t)
		\end{equation}
		has infinitely many solutions $t \in \bN$. Then
		\[
		\nu(\cW(\psi)) =
		\begin{cases}
			0, &\text{if } \displaystyle
			\sum_{t=1}^\infty
			\psi(t) < \infty \\ \\
			1, &\text{if } \displaystyle
			\sum_{t=1}^\infty
			\psi(t) = \infty.
		\end{cases}
		\]
	\end{thm}
	
	We write $\dimH(\cA)$ for the Hausdorff dimension of a Borel set $\cA$, see \cite{Fal2014}. We will make repeated use of the mass distribution principle \cite[Chapter 4]{Fal2014}, so we state a special case of it here for the reader's convenience.
	
	\begin{lemma} [Mass distribution principle]
		\label{MDP}
		Let $F$ be a Borel subset of $\bR^n$, for some $n \in \bN$, and suppose $\mu$ is a Borel probability measure supported on $F$. Assume that, for some $\eps >0 $ and some $s > 0$, we have
		\[
		\mu(U) \ll \diam(U)^s
		\]
		for all Borel sets $U$ such that $\diam(U) \le \eps$. Then
		\[
		\dim_{\mathrm{H}}(F) \ge s.
		\]
	\end{lemma}
	
	Combining Theorems \ref{MainThm} and \ref{LacApprox}, and applying the mass distribution principle,
	delivers the following assertion. 
	
	\begin{thm} \label{LacThm}
		Let $n_1, n_2, \ldots$ be a lacunary sequence of positive integers,
		and let $\psi: \bN \to [0,1]$ with
		\[
		\sum_{t=1}^\infty \psi(t) = \infty.
		\]
		Then there exists $\cG \subseteq \cB$ with $\dimH(\cG) = 2$ such that if $(\bet, \del) \in \cG$ then the inequality
		\eqref{LacIneq} has infinitely many solutions $t \in \bN$. 
	\end{thm}
	
	\noindent
	This is a sharp version of the following result \cite[Theorem 1.11]{CT}.
	
	\begin{thm} [Chow--Technau, 2023+]
		Let $n_1, n_2, \ldots$ be a lacunary sequence of positive integers. Then there exists $\cG \subseteq \cB$ for which $\dimH(\cG) = 2$ such that if $(\bet, \del) \in \cG$ then there exist infinitely many $t \in \bN$ such that
		\[
		\| n_t \bet - \del \| <
		\frac{(\log t)^{3 + \eps}}{t}.
		\]
	\end{thm}
	
	Theorem~\ref{LacThm} follows directly from a more precise counting result, namely Theorem~\ref{CountingThm} below. In order to state it, we first introduce some notation. For any $(\alp,\gam) \in \cB$ and any $N\in\bN$, define
	\begin{equation} \label{def_cN}
		\cN(\alp,\gam,\psi,N)
		= \# \left\{
		t \ls N :  \| n_t \alp - \gam \| 
		\le \psi(n_t)
		\right\}.
	\end{equation}
	Given $\psi: \bN \to [0,1]$, we write
	\begin{equation} \label{PsiDef}
		\Psi(N) = \sum_{t \le N} \psi(n_t).
	\end{equation}
	We denote by $\cR(\psi)$ the set of $(\alp,\gam) \in \bR^2$ for which
	\begin{equation} \label{countlacresult}
		\cN(\alp,\gam,\psi,N) = 2\Psi(N) + O(\Psi(N)^{2/3}
		\left(\log (2+\Psi(N))\right)^3).
	\end{equation}
	The implied constant may depend on $(\alp,\gam)$, but not on $N$. With the above notation at hand, we now state the aforementioned refinement of Theorem~\ref{LacThm}.
	\begin{thm} \label{CountingThm}
		Let $n_1, n_2, \ldots$ be a lacunary sequence of positive integers, and let $\psi: \bN \to [0,1]$. %be such that
		%\begin{equation} \label{divergence}
		%\sum_{t=1}^\infty \psi(n_t) = \infty.
		%\end{equation}
		Then, $\dim_{\mathrm{H}}
		(\cR(\psi) \cap \cB)
		= 2$.
	\end{thm}
	
	Theorem \ref{CountingThm}, in turn, follows from Theorem \ref{MainThm}, the mass distribution principle, and the following result.
	
	\begin{thm}
		\label{GeneralCountingThm} \label{thm_lacunary}
		Fix $A > 2$, and let $\nu$ be a probability measure supported on $[0,1]^2$ satisfying \eqref{LogDecay}. Let $n_1, n_2, \ldots$ be a lacunary sequence of positive integers, and let $\psi: \bN \to [0,1]$.
		Then,
		\[
		\nu(\cR(\psi)) = 1.
		\]
	\end{thm}
	
	\subsection{Multiplicative approximation}
	
	A famous conjecture of Littlewood, from around 1930, asserts that if $\alp, \bet \in \bR$ then
	\[
	\inf_{n \in \bN}
	n \| n \alp \| \cdot \| n \bet \|
	= 0.
	\]
	In pioneering work, Pollington and Velani \cite{PV2000} used Kaufman's measures to demonstrate the following beautiful result.
	
	\begin{thm} 
		[Pollington--Velani, 2000]
		Let $\alp \in \Bad$. Then there exists $\cG \subseteq \Bad$ with $\dimH(\cG) = 1$ such that if $\bet \in \cG$ then
		\[
		\liminf_{n \to \infty} n(\log n)
		\| n \alp \| \cdot \| n \bet \| 
		\le 1.
		\]
	\end{thm}
	
	The assumption that $\alp \in \Bad$ was relaxed to $\alp \in \cK$ in \cite{CZ2021}, where
	\begin{equation} \label{Kdef}
		\cK = \{ \alp \in \bR:
		\sup \{ k^{-1} \log q_k(\alp) \}
		< \infty \}.
	\end{equation}
	Here $q_1(\alp), q_2(\alp), \ldots$ is the sequence of continued fraction denominators for $\alp \in \bR \setminus \bQ$. The set $\cK$ contains $\Bad$ as well as almost all real numbers. The authors of \cite{CZ2021} were also able to make the result completely inhomogeneous, and we state a consequence of \cite[Theorem 1.1]{CZ2021} below.
	
	\begin{thm} 
		[Chow--Zafeiropoulos, 2021]
		Let $\alp \in \cK$ and
		$\gam, \del \in \bR$. Then there exists $\cG \subseteq \Bad$ with
		$\dimH(\cG) = 1$ such that if 
		$\bet \in \cG$ then
		\begin{equation}
			\label{LogThreshold}
			\liminf_{n \to \infty}
			n (\log n) \| n \alp - \gam \|
			\cdot \| n \bet - \del \| \le 1.
		\end{equation}
	\end{thm}
	
	Those authors were able to infer a variant of this in which 
	$(\bet, \del)$ lies in a full-dimensional subset of $\cB$. Their result, which is based on a discrepancy estimate, was recently superseded in 
	\cite[Corollary 1.6]{CT} via
	dispersion.
	
	\begin{thm} [Chow--Technau, 2023+]
		\label{CTlift}
		Let $\alp \in \cK$ and $\gam \in \bR$, and let $\eps > 0$. Then there exists $\cG \subseteq \cB$ with $\dimH(\cG) = 2$ such that if $\bet, \del \in \cG$ then
		\[
		n \| n \alp - \gam \| \cdot
		\| n \bet - \del \| <
		\frac{(\log \log n)^{3+\eps}}
		{\log n}
		\]
		has infinitely many solutions
		$t \in \bN$.
	\end{thm}
	
	Our contribution to the multiplicative theory is to eliminate the 
	$(\log \log n)^{3 + \eps}$ factor.
	
	\begin{thm} \label{MultThm}
		Let $\alp \in \cK$ and $\gam \in \bR$. Then there exists $\cG \subseteq \cB$ with $\dimH(\cG) = 2$ such that if $(\bet, \del) \in \cG$ then
		\eqref{LogThreshold}
		holds.
	\end{thm}
	
	\noindent It is unlikely that this could ever be achieved through dispersion, as discussed in the closing remarks of \cite{CT}. 
	
	As also discussed in those remarks, the proof of Theorem \ref{CTlift} delivers a counting result, that the number of good approximations up to height $N$ is at least a constant times $\log \log N$. We present such a refinement of Theorem \ref{MultThmQ} below. For any $\bet, \del \in \bR$ and any $N \in \bN$, write $\widetilde{\cN}(\alpha,\beta,\gamma,\delta,N)$
	for the number of positive integers $n \le N$ such that
	\begin{equation} \label{strong_Littlewood_ie}
		\| n \alp - \gam \| \cdot
		\| n \bet - \del \| \leq
		\frac{1}
		{n\log n}.
	\end{equation}
	
	Further, write $\cR(\alp, \gam)$ for the set of $(\bet, \del) \in \cG$ such that
	\[
	\widetilde{\cN}(\alpha,\beta,\gamma,\delta,N)
	\gg \log\log N
	\]
	holds for any sufficiently large $N$.
	
	\begin{thm} \label{MultThmQ}
		Let $\alp \in \cK$ and $\gam \in \bR$. Then
		\[
		\dim_H\cR(\alpha,\gamma)
		=2.
		\]
	\end{thm}
	
	Theorem \ref{MultThmQ} follows from Theorem \ref{MainThm}, the mass distribution principle, and the following result.
	
	\begin{thm}
		\label{GeneralMultThmQ}
		Fix $A > 2$, and let $\nu$ be a probability measure supported on $[0,1]^2$ satisfying \eqref{LogDecay}. Let $\alp \in \cK$ and $\gam \in \bR$. Then
		\[
		\nu(\cR(\alpha,\gamma))
		= 1.
		\]
	\end{thm}
	
	Theorems \ref{GeneralCountingThm} and \ref{GeneralMultThmQ} can be applied with $\nu$ as two-dimensional Lebesgue measure. However, the results obtained in this way were already established in \cite{PVZZ2022}.
	
	\bigskip
	
	The aforementioned conjecture about $\Bad_\gam$ remains open. Resolving it would have several desirable consequences. Most of these are similar in flavour to the results presented here. In addition, we would find that $\Bad_\gam$ contains a one-dimensional set of normal numbers.
	
	\subsection{Some consequences for Fourier dimension}
	
	For any Borel set $\cA \subseteq \bR^2$, we write $\cM(\cA)$ for the set of Borel probability measures $\mu$ supported on a compact subset of $\cA$. 
	The \emph{Fourier dimension} of $\cA$ is 
	\[
	\dim_{\mathrm{F}}(\cA) = \sup \{
	s \le 2:
	\exists \mu \in \cM(\cA)
	\quad
	\widehat \mu(\xi) \ll_s 
	(1+|\bxi|)^{-s/2} \:
	(\bxi \in \bR^2) \}.
	\]
	We have the following direct consequence of Theorem \ref{GeneralCountingThm}. 
	
	\begin{cor} Let $n_1, n_2, \ldots$ be a lacunary sequence of positive integers. Then 
		\[
		\dim_{\mathrm{F}}
		([0,1]^2 \setminus \cR(\psi)) = 0.
		\]
	\end{cor}
	
	\noindent This shows that the exceptional set is small in this particular sense. However, we do not know whether $
	\dim_{\mathrm{H}}
	([0,1]^2 \setminus \cR(\psi)) < 1.$
	
	We also have the following direct consequence of Theorem \ref{MainThm}.
	
	\begin{cor} Let $\eta > 0$ be as in Theorem \ref{MainThm}. Then
		\[
		\dim_{\mathrm{F}}(\cB)
		\ge 2\eta.
		\]
	\end{cor}
	
	In view of Remark \ref{rem_concrete_eta}, we have $\dim_{\mathrm{F}}(\cB) \ge 1/15$. We have not tried to fully optimise this. We know from Theorem \ref{MainThm} and the mass distribution principle that $\dim_{\mathrm{H}}(\cB) = 2$. However, we do not know whether $\cB$ is a Salem set, i.e. whether
	\[
	\dim_{\mathrm{F}}(\cB)
	= \dim_{\mathrm{H}}(\cB).
	\]
	This is likely to be a difficult question to answer, since we do not know the answer to the seemingly simpler question of whether $\Bad$ is a Salem set.
	
	\subsection{The new idea}
	
	Our main novelty is in the construction of the measure, that is, in the proof of Theorem \ref{MainThm}. The rest of the paper is similar to \cite{PVZZ2022}. 
	
	In \cite[\S 4]{BHV}, the Ostrowski expansion is used to estimate inhomogeneous approximation errors. Using those findings, we are able to show that $\alp \in \Bad_\gam$ if we impose $R$-periodic constraints on the Ostrowski digits, for some fixed parameter $R \in \bN$. This gives rise to a subset $\cB_{M,R}$ of $\cB$ within which we support $\nu$. If we write down the continued fraction coefficients of $\alp$ and the Ostrowski digits of $\gam$, and left-shift the two sequences $R$ times, then we obtain another element of $\cB_{M,R}$. This symmetry underpins our construction. The remainder of the construction is similar to Kaufman's.
	
	\subsection*{Organisation}
	
	Figure \ref{Leitfaden} illustrates some of the relationships between the new theorems presented above. In Section \ref{prelims}, we provide background knowledge on continued fractions, the Ostrowski expansion, and oscillatory integrals. In Section \ref{SymbolicDynamics}, we present the key new symmetry that underpins our measure construction, and estimate some exponents relating to Hausdorff dimension. Then, in Section \ref{MeasureConstruction}, we construct our new Kaufman-type measures, and show that they have large Frostman dimension. In Section \ref{section_FD}, we show that these measures exhibit polynomial Fourier decay, thereby completing the proof of Theorem \ref{MainThm}. In Section \ref{LacSection}, we establish a refinement of Theorem \ref{GeneralCountingThm} on lacunary approximation. Finally, in Section \ref{MultSection}, we deduce Theorem \ref{GeneralMultThmQ} on multiplicative approximation. As explained earlier, and as depicted in Figure \ref{Leitfaden}, our other results follow from Theorems \ref{MainThm}, \ref{GeneralCountingThm}, and \ref{GeneralMultThmQ}.
	
	\tikzstyle{block} = [rectangle, draw, text centered, rounded corners, minimum height=2em]
	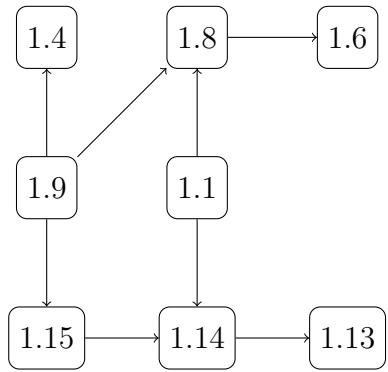
\begin{figure}[!htb] 
		\centering
		\begin{tikzpicture}
			[node distance = 2cm, auto, every text node part/.style={align=center}]
			\node [block] (GenLacCount) {\ref{GeneralCountingThm}};
			\node [block, above of=GenLacCount] (GenLac) {\ref{LacApprox}};
			\draw[->](GenLacCount)--(GenLac);
			\node [block, right of=GenLac] (LacCount) {\ref{CountingThm}};
			\draw[->](GenLacCount)--(LacCount);
			\node [block, right of=LacCount] (Lac) {\ref{LacThm}};
			\draw[->](LacCount)--(Lac);
			\node [block, right of=GenLacCount] (Main) {\ref{MainThm}};
			\draw[->](Main)--(LacCount);
			\node [block, below of=GenLacCount] (GenMultCount) {\ref{GeneralMultThmQ}};
			\draw[->](GenLacCount)--(GenMultCount);
			\node [block, right of=GenMultCount] (MultCount) {\ref{MultThmQ}};
			\draw[->](GenMultCount)--(MultCount);
			\draw[->](Main)--(MultCount);
			\node [block, right of=MultCount] (Mult) {\ref{MultThm}};
			\draw[->](MultCount)--(Mult);
		\end{tikzpicture}
		\caption{Leitfaden}
		\label{Leitfaden}
	\end{figure}
	
	\subsection*{Notation}
	Throughout the text we use the standard $O(\cdot)$ and $o(\cdot)$ notations. Given two complex quantities $A,B$ we shall write $A \ll B$ when  $|A| \ls c |B|$ for some constant $c>0$. We use bold letters such as $\ba, \bb$ to denote tuples of positive integers. For $M \in \bN$, we write $[M] =
	\{1,2,\ldots,M\}$ for the set of the positive integers up to $M$ and $[M]_0 = \{0,1,\ldots, M\}$ for the set of non-negative integers up to $M.$ We write
	$\log^+(x) = \log \max \{x, e \}$.
	
	\subsection*{Funding and acknowledgements} 
	
	SC was supported by EPSRC Fellowship Grant EP/S00226X/2, and by the Swedish Research Council under grant no. 2016-06596. AZ was supported by a postdoctoral fellowship funded by Grant 275113 of the Research Council
	of Norway and also by European Research Council (ERC) under the European Union’s Horizon 2020 Research and Innovation Program, Grant agreement no. 754475. SC thanks Andy Pollington for insights into the Ostrowski expansion, and Thomas Jordan for his enthusiasm. We thank Sanju Velani for his continued encouragement.
	
	\subsection*{Rights} 
	For the purpose of open access, the authors have applied a Creative Commons Attribution (CC-BY) licence to any Author Accepted Manuscript version arising from this submission.
	
	\section{Preliminaries}
	\label{prelims}
	
	As a first preliminary step towards the proofs of the main theorems, we provide the theoretical background for the definition of the measure $\mu$ supported on the set $\cB.$ The desired decay rate of the Fourier transform of $\mu$ will be achieved once we exhibit a certain symmetric structure within $\cB.$ This symmetry will in turn be described in terms of the continued fraction and Ostrowski expansions. 
	
	\subsection{Continued fractions} 
	Given $\alpha\in [0,1) \setminus \bQ,$ we denote its continued fraction expansion by 
	\[
	\alpha = [a_1, a_2, \ldots] = \cfrac{1}{a_1 + \cfrac{1}{a_2 + \dots }}, 
	\]
	where the {\it partial quotients} $a_1, a_2, \ldots$ are positive integers. Given any $k \in \bN$, the $k^{\mathrm{th}}$ \emph{convergent} of $\alpha$ is the rational number
	\[ 
	\frac{p_k}{q_k} = [a_1,\ldots, a_k] = \cfrac{1}{a_1 + \cfrac{1}{ \dots + \cfrac{1}{a_k}}} \, .
	\]
	Here the integers $p_k$ and $q_k$ are relatively prime. The sequences $(p_k)_{k\in\bN}$ and $(q_k)_{k\in\bN}$ satisfy the recursive relations \begin{equation} \label{recursive} 
		p_{k+1} =
		a_{k+1}p_k + p_{k-1},
		\quad
		q_{k+1} =
		a_{k+1}q_k + q_{k-1}
		\qquad (k \ge 0),
	\end{equation}
	where 
	\[
	p_{-1}=1, 
	\qquad
	p_0 = 0, 
	\qquad
	q_{-1}=0,
	\qquad
	q_0=1.
	\]
	It is known that   $p_k/q_k \to \alp$ as $k\to \infty,$ and moreover the convergents of $\alp$ satisfy the inequalities 
	\begin{equation} \label{rate}
		\frac{1}{2q_k q_{k+1}} < \left| \alp - \frac{p_k}{q_k}\right| < \frac{1}{q_k q_{k+1}} \qquad (k\gs 1), 
	\end{equation}
	see  for example  \cite[Chapter 7.5]{NZM1991}. 
	
	Using \eqref{rate}, one can see that the size of the partial quotients is connected with the approximation properties of $\alp.$ In particular, if $a_k \ls M$ for all $k\gs 1,$ then 
	\begin{equation} \label{bad_lower_bound}
		\| q\alp\| \gs \frac{1}{2(M+1)q} \qquad (q\gs 1). 
	\end{equation}
	Indeed, suppose $q, k \in \bN$ with $q_k \ls q < q_{k+1}.$ Then $\|q\alp\| \gs \|q_k\alp\|$, see \cite[Chapter II.2]{RS1992}. By \eqref{rate}, we have
	\[ 
	\|q\alp\| \gs \|q_k \alp\| \gs \frac{1}{2q_{k+1}}.
	\]
	Now $q_{k+1} = a_{k+1}q_{k}+q_{k-1} \ls (M+1)q,$ whence \eqref{bad_lower_bound} follows. 
	
	From the relations \eqref{recursive}, it follows that
	\begin{equation} \label{pkqk}
		p_{k-1}q_k - p_k q_{k-1} = (-1)^{k} \qquad (k\gs 1),
	\end{equation}
	and by induction that
	\begin{equation} \label{cf_inverse}
		\frac{q_{k-1}}{q_k}=[a_k,\dots,a_1] \qquad (k\gs 1).
	\end{equation}
	Regarding the size of the denominators $q_k \,\,  (k\gs 1)$, induction shows that 
	\begin{equation} \label{q_k_size}
		2^{(k-2)/2} \ls q_k \ls (a_1 + 1)\cdots (a_k +1)  \qquad (k\gs 1).    
	\end{equation}
	The continued fraction expansion of $\alpha$ further gives rise to the sequence of {\it remainders} 
	\begin{equation} \label{rk_def}
		r_k = [a_{k+1}, a_{k+2}, \dots] \qquad (k\gs 0).  
	\end{equation} 
	The irrational $\alp$ is linked with the sequences  $(p_k)_{k\in\bN},$ $(q_k)_{k\in\bN}$ and $(r_k)_{k\in\bN}$ by the relation
	\begin{equation} \label{relation_remainders}
		\alp = \frac{p_k + r_{k}p_{k-1}}{q_k + r_{k}q_{k-1}} \qquad (k\geq 1).
	\end{equation}
	As seen from \eqref{recursive}, the denominator $q_k$ of the $k^{\mathrm{th}}$ convergent of $\alp$ is a function of the first $k$ partial quotients of $\alpha,$ and we shall often stress this dependence by writing 
	\[ 
	q_k = q_k(a_1,\ldots,a_k).
	\]
	For future reference, we note that the recursive relations \eqref{recursive} imply the symmetry
	\begin{equation} \label{q_k_inverse}
		q_k(a_1,\ldots,a_k)=q_k(a_k,\ldots,a_1),
	\end{equation}
	which is \cite[Equation (2)]{QR2003}.
	
	Later in the text we will need to use the following property of denominators as functions of the partial quotients \cite[Equation (6)]{QR2003}: for any integers $k,\ell \gs 1$ and for any $a_1,\ldots, a_{k+\ell} \gs 1$, we have 
	\begin{equation} \label{concatenation}
		1 \ls \frac{q_{k+\ell}(a_1, a_2, \ldots, a_{k+\ell})}{q_{k}(a_1,\ldots,a_k) q_{\ell}(a_{k+1},\ldots,a_{k+\ell})} \ls 2.
	\end{equation}  
	Moreover, we define the approximation errors
	\[
	D_k = D_k(\alp) = q_k \alp - p_k
	\qquad (k \ge 0).
	\]
	Regarding the size of the $k^{\mathrm{th}}$ approximation error, the inequalities \eqref{rate} imply that  
	\begin{equation} \label{Dk_size}
		|D_k| \ls   \frac{1}{q_{k+1}}  \qquad (k\gs 1).  
	\end{equation}
	Finally, we note the identity
	\begin{equation} \label{EndpointIdentity}
		\sum_{i=0}^\infty a_{k+2i} |D_{k+2i-1}| = |D_k| \qquad (k \ge 0),
	\end{equation}
	which is \cite[Equation (3.7)]{BHV}.
	
	\subsection{The Ostrowski expansion}\label{Ostrowski_section} Given an irrational $\alp\in [0,1),$ any $\gam \in [-\alp,1-\alp)$ can be  uniquely written in the form 
	\begin{equation} \label{ostrowski}
		\gam = \sum_{k=0}^\infty b_{k+1} D_k, 
	\end{equation}
	where the integers $b_{k+1} = b_{k+1}(\alp,\gamma), \, k=0,1,\ldots$ satisfy
	\begin{align*}
		& 0 \le b_1 \le a_1 - 1, \\
		& 0 \le b_{k+1} \le a_{k+1} \quad (k \in \bN), \qquad \text{ and} \\
		& b_{k+1}  = a_{k+1} \Longrightarrow b_k = 0 \quad (k \in \bN).
	\end{align*}
	The expression in \eqref{ostrowski} is called the \emph{Ostrowski expansion} of $\gamma$ with respect to $\alp$, see \cite{BHV, RS1992}. 
	
	For a fixed irrational $\alp\in [0,1),$ any $n\in \bN$ can be uniquely written in the form  
	\begin{equation} \label{Ostrowski_expansion_integers}
		n= \sum_{k=0}^\infty c_{k+1} q_k, 
	\end{equation}
	where $(c_{k+1})_{k\gs 0}$ is a sequence of integers such that 
	\begin{align*}
		& 0 \le c_1 < a_1,  \\
		& 0 \le c_{k+1} \le a_{k+1}
		\qquad (k \in \bN), \qquad \text{ and} \\
		&  c_{k+1} = a_{k+1} \Longrightarrow c_k = 0 \qquad (k \in \bN).
	\end{align*}
	The right-hand side of \eqref{Ostrowski_expansion_integers} is called \emph{Ostrowski expansion} of the positive integer  $n \in \bN$ with respect to the $\alp\in [0,1) \setminus \bQ.$ Note that the sum in \eqref{Ostrowski_expansion_integers} is actually finite; if we write $K$ for the positive integer such that $q_K \ls n < q_{K+1},$ then \eqref{Ostrowski_expansion_integers} becomes
	\[
	n = \sum_{k=0}^K c_{k+1}q_k.
	\]
	
	For a fixed value of $\alp,$ combining the Ostrowski expansion of the integer $n\gs 1$ with the Ostrowski expansion of the irrational $\gamma$ can provide estimates for the size of the quantity $\|n\alpha -\gam\|,$ as long as $\gamma\notin \bZ\alp + \bZ.$ This condition simply guarantees that $n\alpha -\gamma\notin\bZ$ for any $n\in\bN,$ and in particular $\|n\alp-\gam\|>0$ for any $n\in\bN.$ The following result is \cite[Lemma 4.3]{BHV}; we reproduce it here for the reader's convenience. 
	
	\begin{lemma}\label{nag} Let $\alp \in [0,1)\setminus \bQ,$ let $\gam\in[-\alp,1-\alp)$, and suppose $\gam\notin \alp\bZ + \bZ.$ If $n\in\bN$ then, with reference to the expansions \eqref{ostrowski} and \eqref{Ostrowski_expansion_integers}, put
		\begin{equation} \label{delta_k_def}
			\del_{k+1} = c_{k+1} - b_{k+1} \qquad (k \gs 0).
		\end{equation}
		Let $m\gs 0$ be the smallest integer such that $\delta_{m+1}\neq 0$ and set 
		\[ S = S(n,\alp,\gam) = \sum_{k=m}^{\infty}\delta_{k+1}D_k.  \]
		Then 
		\[
		\|n\alp -\gam\| = \|S\| 
		= \min\{|S|, 1-|S|\}, \qquad 
		|S| = \mathrm{sgn}(\delta_{m+1}D_m)S.
		\]
	\end{lemma}
	
	\begin{cor} \label{cor_nag}
		Under the assumptions of Lemma~\ref{nag}, let us assume in addition that there exists a constant $B>0$ such that $|\delta_k|\leq B$ for all $k\in\bN$. Then
		\[
		\|n\alp-\gam\|\leq \frac{4 B}{q_m}.
		\]
	\end{cor}
	
	\begin{proof} We will use the inequalities
		\begin{equation} \label{useful}
			q_{m+k} \ge F_k q_m \quad (m, k \ge 0),
			\qquad \sum_{k=1}^\infty F_k^{-1} < 4,
		\end{equation}
		where $F_k$ denotes the $k^{\mathrm{th}}$ Fibonacci number.
		The latter is a standard fact. For the former, observe from~\eqref{Dk_size} that $|D_k|\leq q_{k+1}^{-1}$ for all $k\in\bN$. Furthermore, it follows from~\eqref{recursive} that $q_{k+2}\geq q_{k+1}+q_k$ for all $k \ge 0$, so we have $q_{m+k} \ge F_k q_m$.  
		
		Together with Lemma~\ref{nag}, the inequalities \eqref{useful} lead to the upper bound
		\[ \begin{aligned}
			\|n\alp-\gam\| & = \left\|\sum_{k=m}^{\infty}\delta_{k+1}D_k\right\|\leq B\sum_{k=m}^{\infty}|D_k|\\
			&\leq B\sum_{k=m}^{\infty}q_{k+1}^{-1}\leq Bq_m^{-1}\sum_{k=1}^{\infty}\frac{1}{F_k}
			< 4Bq_m^{-1}.
		\end{aligned}
		\]
	\end{proof}
	
	We now present some auxiliary relations that will be used in a later section, when we estimate the size of the Fourier transform of the measure to be constructed. There, given $\alp \in [0,1)\setminus \bQ $ with continued fraction expansion $\alp = [a_1, a_2, \ldots],$ a positive integer $J$ will be fixed and we will need to relate the continued fraction data of $\alp$ with to that of the $J^{\mathrm{th}}$ tail 
	\[ t = 
	[a_{J+1}, a_{J+2}, \ldots]. 
	\] 
	Observe that we can rewrite \eqref{relation_remainders} as
	\begin{equation}
		\alp = \frac{p_J + t p_{J-1}} {q_J + tq_{J-1}}\label{alphaoft}. 
	\end{equation}
	
	\begin{lemma}
		With $\alpha, t$ as above, for any $k\ge 0$ we have 
		\begin{equation}
			\label{DkRelation}
			D_k(t) = - \: \frac{D_{J+k}(\alpha)}{D_{J-1}(\alpha)}.
		\end{equation}
	\end{lemma}
	\begin{proof}
We first observe that   
\begin{equation} \label{Dk_product}
D_k(\alp) = (-1)^k r_0(\alp) r_1(\alp) \cdots r_k(\alp) \qquad (k\gs 0)
\end{equation}
where $r_k$ are the remainders defined in \eqref{rk_def}. This can be shown via induction. For $k \ge 0$, we abbreviate $D_k = D_k(\alp)$. For $k=0$, we have $D_0 = \alp = r_0.$ Assume that for some $k\gs 1$ we have 
\[ 
D_i = (-1)^i r_0 r_1 \cdots r_i \qquad (i\ls k).   
\] 
Then 
\begin{align*}
D_{k+1} 
&= q_{k+1}\alp - p_{k+1} = a_{k+1}D_{k} + D_{k-1}
\\&= (-1)^ka_{k+1}r_0r_1\cdots r_k 
+ (-1)^{k-1}r_0\cdots r_{k-1} 
\\&= (-1)^{k+1} r_0\ldots r_{k-1}
(1- a_{k+1}r_k) 
\\ &= (-1)^{k+1} 
r_0\cdots r_{k-1}r_k(r_k^{-1}
- a_{k+1}) 
\\ &= (-1)^{k+1}
r_0\cdots r_k r_{k+1}.
\end{align*} 

The definition of $t$ implies that $r_k(t) = r_{J+k}(\alp)$ for all $k\gs 0.$  
Together with \eqref{Dk_product}, this gives the desired relation  \eqref{DkRelation}.
\end{proof}

\subsection{Oscillatory integrals} In a later section in the text, we shall make use of the following  lemmas on oscillatory integrals. Note that Lemmas~\ref{lem_ub_nozero_fv} and~\ref{ExpInt} are similar to  \cite[Lemmas~5.1 and~5.2]{QR2003},  the principal difference being that we integrate over arbitrary subintervals of $[0,1]$. One might notice that our---a priori more general---Lemmas~\ref{lem_ub_nozero_fv} and~\ref{ExpInt}   can be deduced from \cite[Lemma~5.1]{QR2003} and~\cite[Lemma~5.2]{QR2003} respectively via a linear change of variables. We nonetheless provide a full proof, as  the details are omitted in \cite{QR2003}.

\begin{lemma} \label{lem_ub_nozero_fv}
Let $a,b > 0$, and suppose $F \in C^2([-1,1])$ satisfies 
\[
|F'(t)| \ge a, \quad |F''(t)| \le b \qquad (-1 \le t \le 1).
\]
Then, for any subinterval $\cI$ of $[-1,1]$, we have
\[
\int_\cI e(F(t)) \d t \ll \frac1a + \frac{b}{a^2}.
\]
\end{lemma}

\begin{proof}
We denote by $\alpha$ and $\beta$ the endpoints of $\cI$, so that $\cI=[\alpha,\beta]$. Integration by parts yields
\[
\begin{aligned}
\int_\cI e(F(t)) \d t &=\int_\cI e(F(t))F'(t)\frac{1}{F'(t)} 
\d t=\int_\cI 
\left(\frac{\d}{\d t} e(F(t)) \right)
\frac{1}{2\pi i F'(t)} \d t\\
&=\left.\frac{e(F(t))}{2\pi i F'(t)}\right|_{\alpha}^{\beta}
+ \int_\cI e(F(t))\frac{F''(t)}{2\pi i \left(F'(t)\right)^2} \d t.
\end{aligned}
\]
Thus, using the triangle inequality, we find that
\[
\int_\cI e(F(t)) \d t \ll \frac{1}{\min_{t\in\cI}\left|F'(t)\right|}
+
\frac{\max_{t\in\cI}\left|F''(t)\right|}{ \min_{t\in\cI}\left|F'(t)\right|^2}\leq \frac1a + \frac{b}{a^2}.
\]
\end{proof}

\begin{lemma} \label{lem_ub_nozero}
Let $k \in \bN$ and $\lam > 0$. Let $\cI \subseteq \bR$ be an open interval and suppose $F \in C^k(\cI)$ satisfies $|F^{(k)}(t)| \ge \lam$ for all $t\in \cI$. In the case $k = 1$, assume further that $F'$ is monotonic. Then
\[
\int_\cI e(F(t)) \d t \ll_k \lam^{-1/k}.
\]
\end{lemma}

\begin{proof}
This follows from \cite[Chapter VIII, \S 1, Proposition 2]{Ste1993}.
\end{proof}

\begin{lemma} \label{ExpInt} Let $m\in\mathbb{N}$,  $f\in C^2([0,1])$ and assume that $f''$ has at most $m$ zeros in $[0,1]$. Let  $a,b\in\mathbb{R}$ be real numbers such that
\begin{equation} \label{a_less_than_b}
0\leq a\leq mb.
\end{equation}
Suppose
\[
f'(t) = (A t + B) G(t),
\]
where $A, B \in \bR$ with $A \ne 0$  and 
\[
|G(t)| \ge a, \quad |G'(t)| \le b \qquad (0 \le t \le 1).
\]
Then, for $0 \le X \le Y \le 1$, we have 
\begin{equation} \label{lem_ub_zero_result}
\int_X^Y e(f(t)) \d t \ll \frac{m b}{a^{3/2} |A|^{1/2}}.
\end{equation}
\end{lemma}

\begin{proof}
First of all, note that the left-hand side of \eqref{lem_ub_zero_result} is trivially  bounded above by $1$. 
Because of this, in the rest of the proof it is enough to consider parameters such that the right-hand side of~\eqref{lem_ub_zero_result} is strictly less than 1, that is 
\begin{equation} \label{rhs_is_small}
\frac{b}{|A|^{1/2} a^{3/2}} < \frac{1}{m} \cdot 
\end{equation}
Further, note that the right-hand side of~\eqref{lem_ub_zero_result} is decreasing with respect to the parameter $a$. So, without loss of generality, we can assume that
\begin{equation} \label{a_is_min_g}
a = \min \{ G(t): 0 \le t \le 1 \} = G(x_0)
\end{equation}
for some $x_0 \in [0,1].$ 

Define
\[
F(x) = f(x-B/A)  \qquad (B/A \le x \le 1+B/A).
\]
Then
\[
\int_X^Y e(f(t)) \d t =   \int_{X+B/A}^{Y+B/A} e(F(x)) \d x.
\]
We compute that
\[
F'(x) = A x G(x - B/A), \qquad
F''(x) = A G(x - B/A) + A x G'(x - B/A).
\]
Put $\cI = [X+B/A, Y+B/A]$ and $c = a/(2b)$. For $|x| \le c$, we have
\[
|F''(x)| \ge |A| (a - cb) = |A|a/2,
\]
and Lemma~\ref{lem_ub_nozero} 
yields
\[
\int_{\cI \cap [-c,c]} e(F(x)) \d x \ll (|A| a)^{-1/2}.
\]

Next, decompose $\cI \setminus [-c,c]$ into a disjoint union
\[
\cI \setminus [-c,c] = \cI_1 \cup \cdots \cup \cI_{m+2}
\]
of intervals on which $F'$ is monotonic, and consider some 
$j \in [m+2]$.
For all $x\in \cI \setminus 
[-c, c]$, we have
\[
|F'(x)| > |A|ca \gg |A|a^2/b.
\]
Now Lemma~\ref{lem_ub_nozero} yields
\[
\int_{\cI_j} e(F(x))
\mathrm{d}x \ll \frac{b}{|A|a^2}.
\]
Combining the estimates above and applying~\eqref{rhs_is_small} and \eqref{a_less_than_b}, we conclude that 
\begin{align*}
\int_X^Y e(f(t)) \d t &\ll (|A|a)^{-1/2} + (m+2)\frac{b}{|A|a^{2}} \\ &\leq 
\Big(1 + \frac{m+2}{m}\Big) 
\frac{1}{|A|^{1/2} a^{1/2}} 
\ll \frac{mb}{|A|^{1/2} a^{3/2}}.
\end{align*}
\end{proof}

\section{Symbolic dynamics}
\label{SymbolicDynamics}

In the present section we introduce large subsets of $\cB$, defined via the Ostrowski expansion, within which we will support our measure $\nu$ in Theorem \ref{MainThm}. Their self-similar structure will be the key to establishing the decay property of the Fourier transform of $\nu.$ 

\subsection{A self-similarity}

Let $M \ge 2$ and $R \ge 3$ be integers. We write 
\[
\cF_M = \left\{\alp \in [0,1)\setminus \bQ : 1\ls a_k \ls M \text{ for all }  
k \in \bN \right\}  
\]
for the set of irrational numbers $\alp$ in $[0,1)$ whose partial quotients are all bounded above by $M$. By a standard characterisation of badly approximable numbers, we have
\[ 
\Bad \cap [0,1) = \bigcup_{M=2}^\infty \cF_M. 
\]
For $\alp \in \bR \setminus \bQ$, denote by $\cO_{R}(\alp)$
the set of $\gam \in [-\alp,1-\alp)$ such that the sequence of digits associated with the Ostrowski expansion of $\gam$ with respect to $\alp$ as in \eqref{ostrowski} satisfies
\begin{equation}
\label{PeriodicConstraint} 
1 \le b_k \le a_k-1
\qquad (k \mmod R \in \{0,1\}).
\end{equation}
Finally, define
\[
\cB_{M,R} = \{ (\alp, \gam) \in \bR^2: \alp \in \cF_M \text{ and } \gam \in \cO_{R}(\alp) \}.
\]

\begin{lemma} We have $\cB_{M,R} \subseteq \cB$.
\end{lemma}

\begin{proof} Let $(\alp,\gam) \in \cB_{M,R}$. Let $n\in \bN$, and write
\[
n = \sum_{k=0}^K c_{k+1} q_k
\asymp q_K 
\]
for the Ostrowski expansion of $n$ with respect to $\alp$. As in Lemma~\ref{nag}, we let $m \gs 0$ be minimal such that $\del_{m+1} \ne 0,$ where $\del_k$ is as in \eqref{delta_k_def}. Note from Lemma \ref{nag} that 
\[
\| n \alp - \gam \| = \min \{ |S|, 1-|S| \},
\]
where
\[
S = \sum_{k=m}^\infty \del_{k+1} D_m.
\]
It therefore suffices to find lower bounds for $|S|$ and $1-|S|$. It follows from \cite[Lemma 4.4]{BHV} that there exists some 
$\ell \ls \max\{2, K-m + 1\}$ such that $|S| \gs |D_{m+\ell +1}|.$ Since $(\alp,\gam) \in \cB_{M,R}$ implies that $m \le K + R$, this gives
\[
|S| \ge |D_{m+1+\ell}| \gg 
\min \{ q_m^{-1}, q_K^{-1} \}
\gg q_K^{-1} \gg n^{-1}. 
\]
Regarding $1-|S|,$ it is shown in \cite[Lemma 4.5]{BHV} that for some positive integer $L \le K+2$, 
\[ 
1 - |S| \gs |D_L| \gg q_L^{-1} \gg  n^{-1}.
\]
Combining, we deduce that 
\[
\| n \alp - \gam \| \gg n^{-1},
\]
where the implied constant depends at most on $M$ and $R$, whence  $\cB_{M,R} \subseteq \cB$.
\end{proof}

The \emph{$R$-fold left shift operator} acts on $\cB_{M,R}$ by
\[
\sig^R(\alp, \gam)
= (\alp', \gam')
\]
where, if
$
\alp = [a_1,a_2,\ldots]
$
and $\gam$ has Ostrowski expansion $\displaystyle \sum_{k \ge 0} b_{k+1} D_k$ then
\[
\alp' = [a_{R+1},\ldots],
\qquad
\gam' = \sum_{k \ge 0} 
b_{R+k+1} D_k.
\]
Symbolically, we have
\[
((a_k)_{k=1}^\infty,
(b_k)_{k=1}^\infty)
\mapsto
((a_{R+k})_{k=1}^\infty,
(b_{R+k})_{k=1}^\infty).
\]
This self-similarity is needed for the construction of our measure $\mu$.

\subsection{Exponents relating to Hausdorff dimension} \label{subsection_HD}

Let $m$ be a positive multiple of $R$, and let $\ba \in \bN^m$. 

\begin{defn} \label{Vdef}
For $k \in [m]$, denote by $\cV_k(\ba)$ the set of 
\[
\bb = (b_1,\ldots,b_k) \in \bZ^k
\]
such that
\begin{align} 
&1 \le b_k \le a_k - 1 \qquad 
(k \mmod R \in \{0,1\}),
\notag \\
&0 \le b_k \le a_k \qquad (k \in [t]), \notag
\\ &b_{k+1} = a_{k+1} \implies b_k = 0 \qquad (k \in [t-1]).\notag
\end{align}
Let $\kap_{M,R,m} > 1/2$ be the value of $\kap$ for which
\begin{equation} 
\label{def_kappa_MRm}
\sum_{\ba \in [M]^m} q_m(\ba)^{1-2 \kap} |\cV_m(\ba)| = 1,
\end{equation}
and let $\ome_{M,m} > 0$ be the value of $\ome$ for which
\[
\sum_{\ba \in [M]^m} 
q_m(\ba)^{-2 \ome} = 1.
\]
\end{defn}
Good \cite{Good} showed that
\[
\ome_{M,m} \to \dim_H \cF_M 
\qquad (m \to \infty),
\]
see \cite[Theorem 3.1]{QR2003}. Jarn\'ik \cite{Jar1928} showed that
\[
\dim_H \cF_M \to 1 
\qquad (M \to \infty),
\]
so
\begin{equation} \label{GoodJarnik}
\ome_{M,m} = 1 + o(1) 
\qquad (\min \{ M, m \} \to \infty).
\end{equation}

For $k \in [m]$, write $v_k = |\cV_k(\ba)|$, and note that if $v_m > 0$ and $k \mmod R \in \{0,1\}$ then $a_k \ge 2$. Moreover, we have
\[
v_1 = a_1 - 1, \qquad v_2 = a_2(a_1 - 1),
\]
Further, if $k \ge 3$ and $k \mmod R \in \{0,1\}$ then
\[
v_k = (a_k - 1)v_{k-1}.
\]

\begin{lemma} Let $k \in [3,m]$ be an integer such that $k \mmod R \notin \{0,1\}$. Then
\[
v_k =
\begin{cases}
a_k v_{k-1}, &\text{if }
k \equiv 2 \mmod R \\
a_k v_{k-1} + v_{k-2}, &\text{if } k \not \equiv 2 \mmod R.
\end{cases}
\]
\end{lemma}

\begin{proof} Write
\[
v_{k-1} = y_{k-1} + z_{k-1},
\]
where $y_{k-1}$ counts $\bb \in \cV_{k-1}(\ba)$ such that $b_{k-1} \ne 0$, and $z_{k-1}$ counts $\bb \in \cV_{k-1}(\ba)$ such that $b_{k-1} = 0$. Then
\[
v_k = a_k y_{k-1} + (a_k + 1) z_{k-1}
= a_k v_{k-1} + z_{k-1},
\]
and
\[
z_{k-1} = 
\begin{cases}
0, &\text{if }
k \equiv 2 \mmod R \\
v_{k-2}, &\text{if } k \not \equiv 2 \mmod R.
\end{cases}
\]
\end{proof}

We now see that if $k \in \bN$ then
\[
v_{k+2} \ge v_k,
\]
and hence
\[
v_{2k} \ge 2^{\lfloor k/2 \rfloor}.
\]

\begin{lemma} \label{vq} 
Assume that for $r \mmod R \in \{0,1\}$ we have $a_r \ge 2$. Then
\[
v_k \le q_k(\ba) \le 3^{n_k} v_k \qquad (k \in [m]),
\]
where $n_k$ denotes the number of $j \in [k]$ such that $j \mmod R \in \{0,1,2\}$.
Consequently
\[
v_k \le q_k(\ba) \le 
27^{\lceil k/R \rceil} v_k 
\qquad (1 \le k \le m).
\]
In particular, there exists an absolute constant $c > 0$ such that
\[
v_m \le q_m(\ba) \le v_m^{1+c/R}.
\]
\end{lemma}

\begin{proof} We induct on $k$. We have
\[
q_1 = a_1, \qquad q_2 = a_2 a_1 + 1
\]
and
\[
v_1 = a_1 - 1 \in [q_1/2,q_1],
\qquad
v_2 = a_2 (a_1 - 1) \in [q_2/3,q_2].
\]
Next, let $k \in [3,m] \cap \bZ$, and assume that the result holds with $k$ replaced by $k-1$ or $k-2$. Then
\[
v_k \le a_k v_{k-1} + v_{k-2} \le a_k q_{k-1} + q_{k-2} = q_k.
\]
If $k \mmod R \notin \{0,1,2\}$ then
\[
q_k = a_k q_{k-1} + q_{k-2}
\le 3^{n_k} (a_k v_{k-1} + v_{k-2}) = 3^{n_k} v_k.
\]
If $k \mmod R \in \{0,1,2\}$ then
\[
q_k =  a_k q_{k-1} + q_{k-2} \le (a_k + 1) q_{k-1} \le 3^{n_k - 1} (a_k + 1) v_{k-1} \le 3^{n_k} v_k.
\]
\end{proof}

\begin{lemma} \label{lem_kappa_tends_to_two} We have
\[
\kap_{M,R,m} = 2 + o(1) \qquad
(\min\{M, R, m\} \to \infty).
\]
\end{lemma}

\begin{proof}
It is convenient to abbreviate 
$\kap = \kap_{M,R,m}$ and 
$\ome = \ome_{M,m}$ in the following calculation.
By Lemma \ref{vq}, we have
\begin{align*}
\sum_{\ba \in [M]^m} 
q_m(\ba)^{-2\ome}
&= 1
= \sum_{\ba \in [M]^m} 
q_m(\ba)^{1-2\kap} |\cV_m(\ba)|
\\&= \sum_{\ba \in [M]^m} q_m(\ba)^{2-2\kap+O(1/R)},
\end{align*}
and so
\[
\kap = 1 + \ome + o(1) \qquad (R \to \infty).
\]
Combining this with \eqref{GoodJarnik} completes the proof.
\end{proof}

\section{An analogue of Kaufman's measures}
\label{MeasureConstruction}

\subsection{Construction of the measures} \label{subsection_ConstructionOfMeasures}

We are now in position to construct the probability measures $\nu$ from Theorem \ref{MainThm}. Let 
$\eps \in (0,1/4)$ and 
$\kap \in (1,2)$ be fixed. The parameters $m,M, R$ are as in \S \ref{SymbolicDynamics}: $M\gs 2$ and $R\gs 3$ are integers, and $m\gs R$ is a multiple of $R.$ 
For any such integer $m$, we define 
\[
T_m(\kap) = \sum_{\ba \in [M]^m} q_m(\ba)^{1-2\kap}|\cV_m(\ba)|.
\]
Here we write $|\cV_m(\ba)|$ for the cardinality of the set $\cV_m(\ba)$ introduced in Definition \ref{Vdef}.  
In view of Lemma~\ref{lem_kappa_tends_to_two}, if $M,R,m$ are large enough, then we can ensure that $\kappa_{M,R,m}$ is arbitrarily close to 2, and we choose these parameters large enough so that
\begin{equation} \label{kappa_MRm_bigger_than_kappa}
\kappa_{M,R,m}>\kappa.
\end{equation}
The definition of $\kappa_{M,R,m}$ and the inequality
\eqref{kappa_MRm_bigger_than_kappa} 
imply that
\[
T_m(\kap) \to \infty \qquad (m \to \infty).
\]
We take $m \gs 100$ large enough to ensure that
\begin{equation} \label{T_m_is_large}
T_m(\kap) > 4.
\end{equation}

Consider the set 
\[ 
\cJ = \left\{ (\ba,\bb) \in [M]^m \times [M]_0^m : \bb \in \cV_m(\ba) \right\}. 
\]
The product set $\cJ\times \cJ \times \ldots $ can be naturally identified with the set of pairs $(\alp, \gamma)$ with $0 \le \alp < 1$ and $-\alp \le \gamma < 1-\alp $ such that $\alp \in \cF_M$ and the Ostrowski expansion of $\gamma$ with respect to $\alp$ satisfies 
\eqref{PeriodicConstraint}. 
In particular, we shall not make any distinctions between measures defined on $\cJ \times \cJ \times \ldots$ and on the aforementioned set.
Define a probability measure $\lambda_m$ on $\cJ$ by setting 
\begin{equation} \label{def_lambda_m}
\lam_m(\ba,\bb) = \frac{1}{T_m(\kap)}q_m(\ba)^{1-2\kap} \qquad 
((\ba, \bb) \in \cJ).
\end{equation}
Note that we write $\lam_m(\ba,\bb)$ for notational convenience, instead of the more formal 
$\lam_m(\{ (\ba, \bb)\}).$

Next, we define $\sig_m(\kap)$ so that $m \sig_m(\kap)$ is the mean of the random variable $\log q_m(a_1,\ldots,a_m)$ with respect to the measure $\lam_m$. In other words, 
\[ 
m \sig_m(\kap) = \sum_{(\ba,\bb) \in \cJ} \log q_m(\ba) \lambda_m(\ba,\bb). 
\]
The assumption that $m \geq 100$, along with \eqref{q_k_size}, imply that if $\mathbf{a}\in [M]^m$ then
\[
\log q_m(\mathbf{a})\geq m/3,
\]
whence $\sigma_m(\kap) \geq 1/3$. Henceforth, we assume that $m$ is large enough so that
\begin{equation} \label{m_conditions}   \max \left(\frac{\log 2}{m\sigma_m(\kap)   }, \frac{3}{m} \right) < \varepsilon. 
\end{equation}
In preparation for Section~\ref{section_FD}, we also assume that $R$ is sufficiently large to ensure that
\[
\frac{3 \log 27}{R} (1+c)
< \eps,
\]
where $c$ is as in Lemma \ref{vq},
whence
\begin{equation} \label{R_is_large}
\frac{c}{R} < \eps,
\qquad
\frac{\log 27}{ R\sigma_m(\kappa)} <\varepsilon.
\end{equation}

On the space $\cJ \times \cJ \times \ldots$, we define the probability measure 
\[  
\lambda_m \times \lambda_m \times \ldots,
\]
and let $(Y_n)_{n=1}^{\infty}$ be the sequence of random variables on this probability space defined by
\[
Y_1 = q_m(a_1,\ldots,a_m), \quad Y_2 = q_m(a_{m+1},\ldots,a_{2m}),\ldots.
\]
By the law of large numbers, there exists $j_0 \gs 1$ such that the  set 
\[
\cE = \cE(j_0) \subseteq \cJ^{j_0}
\]
on which
\begin{equation} \label{Edef}
\left|
\frac{\log Y_1 + \cdots + \log Y_{j_0}}{j_0}
- \bE(\log Y_1)
\right| \le \eps \bE(\log Y_1)
\end{equation}
has measure \begin{equation} \label{cEisLarge}
\lam^{j_0}(\cE) > 1/2.
\end{equation} 
We fix such a value of $j_0$, and write 
$\bP = \lam^{j_0}$ as well as
\[
J_0 = j_0m.
\]
Using \eqref{concatenation} and \eqref{Edef}, we deduce that if $N \gs 1$ and $(\ba,\bb) \in \cE^N$ then
\begin{equation} \label{ie_logQJ}
\begin{aligned}
& |\log q_{J_0 N}(a_1,\ldots,a_{J_0 N}) - J_0 N \sig_m(\kap)| \\
& \qquad \le \eps J_0 N \sig_m(\kap) + (j_0N-1) \log 2 \\
& \qquad \le \left(\eps + \log 2\frac{j_0 N -1}{J_0 N \sigma_m(\kappa) } \right) J_0 N \sig_m(\kap)\\
& \qquad < 2\eps J_0 N \sig_m(\kap),
\end{aligned}
\end{equation}
and consequently
\begin{equation} \label{Qest}
Q^{1-2\eps} \le 
q_{J_0N}(a_1,\ldots,a_{J_0N}) \le Q^{1+2\eps},
\end{equation}
where $Q = \exp\left(J_0 N \sig_m(\kap)\right)$.

Let $\lam$ be the probability measure induced by $\bP$ on $\cE$, i.e.
\[
\lambda(A) = \frac{1}{\bP(\cE)}\bP(A\cap \cE) 
\qquad 
(A \subseteq \cJ^{j_0}).
\]
The normalising factor above is the constant
\begin{equation} 
\label{def_c_lambda}
c_{\lambda}:=\frac{1}{\bP(\cE)}, 
\end{equation}
which will appear frequently in subsequent calculations. By~\eqref{cEisLarge}, we have
\[
1\leq c_\lambda < 2.
\]
As $\cE \subseteq \cJ^{j_0}$,
we may also view $\cE \times \cE \times \ldots$ as a set of pairs $(\alp, \gam).$ Finally we define the probability measure $\nu$ by setting
\begin{equation} \label{def_nu}
\nu = \lam\times\lam\times\ldots.
\end{equation}

\subsection{Frostman dimension} 

The main aim of this subsection is to estimate the Frostman dimension~\cite{FFK2021, KPV2022} of the measure $\nu$. Recall that this dimension is defined as the supremum of $s \in \bR$ such that
\[
\nu(B_r(\bx)) \ll r^s
\qquad (0 < r < 1, \: \bx \in \bR^2),
\]
where $B_r(\bx)$ denotes the ball of radius $r$ centred at $\bx$. By the mass distribution principle (Lemma \ref{MDP}), the Frostman dimension of a Borel probability measure is bounded above by the Hausdorff dimension of its support.

In the sequel, it will be convenient to refer to certain subsets of $\bR$ and $\bR^2$ as \emph{cylinder sets}. Given $N\gs 1$ and $\ba = (a_1,\ldots, a_N) \in \bN^N,$ we define the cylinder set $C(\ba)$ as 
\[ 
C(\ba) = \{ \alpha \in [0,1) : a_i(\alp) = a_i
\quad (1\ls i \ls N) \}. 
\]
In other words, the set $C(\ba)$ comprises elements of $[0,1)$ whose continued fraction expansion starts with $\ba$. Further, given $N\gs 1,$ $\ba\in \bN^N$ and $\bb\in\cV_N(\ba)$, we define the cylinder set $C(\ba,\bb)$ as
\[ 
\{ (\alpha,\gamma) \in \bR^2 : \alp \in C(\ba), \: -\alpha
\le \gamma < 1-\alpha, \:
b_i(\alp,\gam) = b_i
\quad (1\ls i \ls N) \}. 
\] 
In other words, the set $C(\ba,\bb)$ comprises $(\alpha,\gamma)$ such that $\alp \in C(\ba),$  $-\alpha \le \gamma
< 1-\alpha$,
and the Ostrowski expansion~\eqref{ostrowski} of the coordinate $\gamma$ with respect to $\alp$  starts with the word $\bb$.
\begin{lemma} \label{lem_size_of_cylinder_set}
We use the notation and assumptions of this section, and assume that $\kappa>1/2$, so that $1-2\kappa<0$. 
Put
$
Q_0 =
\exp(J_0 \sigma_m(\kappa)),
$
and let $N \gs 1$. Then, for any cylinder set $C=C(\ba,\bb)$ with $(\ba,\bb)\in\cE^N$, we have
\[
\frac{1}{2^{j_0-1}} 
\frac{Q_0^{(1-2\kappa)
	(1-2\varepsilon)N}}
{T_m(\kappa)^{j_0 N}} \leq \nu(C)\leq 2^{N}
\frac{Q_0^{(1-2\kappa)(1+2\varepsilon)N}}
{T_m(\kappa)^{j_0 N}}. 
\]
\end{lemma}

\begin{proof}
Writing $\ba = (\ba_1,  \ldots, \ba_N)$ and $\bb = (\bb_1, \ldots, \bb_N)$ where 
\[
(\ba_i,\bb_i) \in \cE 
\qquad (i\ls N),
\]
the definition of $\nu$ in \eqref{def_nu} gives 
\[
\nu(C) = \lambda(\ba_1,\bb_1)
\cdots \lambda(\ba_N, \bb_N).
\]
The definition of $\lambda$  and~\eqref{cEisLarge} imply that   
\begin{equation*} \label{lambda_and_lambda_m0_are_comparable}
\lambda_m^{j_0}(\ba_i,\bb_i) \leq \lambda(\ba_i,\bb_i)\leq 2\lambda_m^{j_0}(\ba_i,\bb_i)  \qquad ( i\ls N).
\end{equation*}
Combining these facts, we infer
\begin{equation} \label{lem_size_of_cylinder_set_db_nu}
\lambda_m^{j_0 }(\ba_1, \bb_1)\cdots   \lambda_m^{j_0}(\ba_N,\bb_N) \leq \nu(C)\leq 2^N  \lambda_m^{j_0 }(\ba_1, \bb_1)\cdots   \lambda_m^{j_0}(\ba_N,\bb_N).
\end{equation}
Further, from the definition~\eqref{def_lambda_m}, we have
\begin{equation*} 
\lambda_m^{j_0 }(\ba_i, \bb_i) = \prod_{k=0}^{j_0-1} \frac{q_m(a_{(i-1)J_0+km+1},\dots,a_{(i-1)J_0+(k+1)m})^{1-2\kap}}{T_m(\kap)}      \quad (i\ls N).
\end{equation*}
Using the inequality \eqref{concatenation}, we deduce that for $i=1, \dots, N$,
\begin{equation}
\begin{aligned}
	\label{lem_size_of_cylinder_set_db_denominator}
	\frac{1}{2^{j_0-1}} \frac{ q_{J_0}(a_{(i-1)J_0+1},\dots,a_{iJ_0})^{1-2\kap}}{T_m(\kap)^{j_0}} &\leq  \lam_m^{j_0}(\ba_i, \bb_i) 
	\\& \leq  \frac{q_{J_0}(a_{(i-1)J_0+1},\dots,a_{iJ_0})^{1-2\kap}}{T_m(\kap)^{j_0}} \cdot 
\end{aligned}
\end{equation}
We further estimate, by using~\eqref{ie_logQJ}, that for $i=1, \dots, N$ we have
\begin{equation} \label{lem_size_of_cylinder_set_db_q_m}
Q_0^{(1-2\varepsilon)(1-2\kap)}
\leq
q_{J_0}(a_{(i-1)J_0+1},\dots,a_{iJ_0})^{1-2\kap}\leq Q_0^{(1+2\varepsilon)(1-2\kap)}.
\end{equation}
The assertion of the lemma follows by putting together~\eqref{lem_size_of_cylinder_set_db_nu}, \eqref{lem_size_of_cylinder_set_db_denominator} and~\eqref{lem_size_of_cylinder_set_db_q_m}.
\end{proof}

The following is a two-dimensional variant of \cite[Proposition 4.1]{QR2003}.  

\begin{lemma} \label{lem_FD}  Let $\nu$ be the measure constructed in Subsection~\ref{subsection_ConstructionOfMeasures}, and assume in addition that the parameter $\kap$ from that section lies in the range $(3/2,2)$. Then we have the following.

\begin{enumerate}[(a)]
\item \label{point_sizeQJ} For all positive integer multiples $J$ of $J_0$, the inequalities
\begin{equation} \label{denominator_bounds}
	Q^{1-2\eps} \le q_J(\alp) \le Q^{1+2\eps},
	\qquad
	\frac{Q^{1-2\eps}}{2M}
	\le q_{J-1}(\alp) \le Q^{1+2\eps} 
\end{equation}
hold for $\nu$-almost all $(\alp,\gam),$ where $Q = \exp(J \sig_m(\kap))$.

\item \label{point_dimFrostman} There exists a constant $c_1>0$ depending only on $M, J_0 \gs 1$ and $\sig_m(\kap)$ such that if $\bx \in \bR^2$ and $r \in (0,1]$ then
\[
\nu(B_r(\bx)) \ls c_1 r^{2\kappa-2-4(2-\kap)\varepsilon}.
\]
\end{enumerate}
\end{lemma}

\begin{proof}
\eqref{point_sizeQJ} The first inequality was already shown in \eqref{Qest}. The second inequality then follows from the bounds
\[
1 \le
\frac{q_J(\alp)}{q_{J-1}(\alp)} \le 2M.
\]

\eqref{point_dimFrostman}  Fix $\bx\in\bR^2$ and $r>0.$ If $\nu(B_r(\bx))=0$ then the claim is straightforward, so in the rest of the proof we assume that $\nu(B_r(\bx))>0$. Define
\begin{equation} \label{def_Q0_N}
Q_0:=\exp(J_0 \sigma_m(\kappa)) \qquad \text{ and  } \qquad N:=\left\lfloor\frac{\log r}{-\log Q_0}\right\rfloor+1.
\end{equation}
These definitions imply  
\begin{equation} \label{r_level}
Q_0^{-N} < r \le Q_0^{1-N}.
\end{equation}
Let $C=C(\ba,\bb)$ be a cylinder set with $(\ba,\bb)\in\cE^N$  such that
\begin{equation} \label{lem_FD_nontrivial_intersection}
\nu\left(B_r(\bx)\cap C\right)>0.
\end{equation}
Setting  $J=J_0 N$ and $Q=\exp(J\sig_m(\kap))=Q_0^N$ as above, point~\eqref{point_sizeQJ} of the lemma gives
\begin{equation} \label{lem_FD_size_q}
Q^{1-2\varepsilon}
\le q_{J}(\ba)\le Q^{1+2\varepsilon}.
\end{equation}

We aim to provide an upper bound for the number of cylinder sets $C(\ba,\bb)$, where $(\ba,\bb) \in \cE^N$, satisfying \eqref{lem_FD_nontrivial_intersection}. We do so by first bounding the number of possible tuples $\ba$ and then the number of possible tuples $\bb$. For each $\ba\in [M]^{J}$, the cylinder set $C(\ba)$ has length 
\[
|C(\ba)| \gs \frac1{2Q^{2(1+2\eps)}},
\]
for the elements of $C(\ba)$ are given by \eqref{alphaoft} for $t \in [0,1)$.
Hence, there exist at most
\begin{equation} \label{lem_FD_ub_cylinders_a}
2r\cdot 2Q^{2(1+2\varepsilon)}+1 \ls 4Q^{1+4\varepsilon} +1 \ls 5Q^{1+4\varepsilon}
\end{equation}
tuples $\ba$ such that the cylinder set $C(\ba,\bb)$ satisfies \eqref{lem_FD_nontrivial_intersection} for some $\bb$.

Now for a fixed $\ba= (a_1, \ldots, a_{J}) \in [M]^{J}$, we find an upper bound for the number of different $J$--tuples $\bb\in \cV_{J}(\ba)$ such that $C(\ba,\bb)$ satisfies \eqref{lem_FD_nontrivial_intersection}. Set $$\alpha=\alpha(\ba) = [a_1,\ldots, a_{J}, M, M, \ldots].$$ 
For $\bb = (b_1, \ldots, b_{J}) \in \cV_{J}(\ba)$,  let 
\[ n = n(\ba,\bb) =  b_{1}q_0(\alpha) + b_2q_1(\alpha) + \ldots + b_{J}q_{J-1}(\alpha) \]
be the integer whose digits in the  Ostrowski expansion~\eqref{Ostrowski_expansion_integers} with respect to $\alpha$ are the entries of $\bb$. By Corollary~\ref{cor_nag},
if $(\alpha,\gamma) \in C(\ba,\bb)$ then
\[
\|\gamma-n(\ba,\bb)\alpha\| < 4M/ q_{J}(\ba). 
\]
This means that $\gamma$ belongs to an interval in $\bR/\bZ$ centred at $n(\ba,\bb)\alpha$ and of length at most $8 M/q_{J}(\ba)$. Note from 
\eqref{useful} that if $\bb_1, \bb_2 \in \cV_{J}(\ba)$ then
\[
|n(\ba,\bb_1)-n(\ba,\bb_2)|
\le 4M q_{J}(\ba),
\]
so by \eqref{bad_lower_bound} the centres of the aforementioned intervals are separated by a distance of at least 
\begin{align*} \|n(\ba,\bb_1)\alpha-n(\ba,\bb_2)\alpha\| &\gs  (2(M+1)|n(\ba,\bb_1)-n(\ba,\bb_2)|)^{-1}\\  
& \ge \left(8M(M+1) q_{J}(\ba)\right)^{-1}.
\end{align*}
Thus, any interval of length $2r$ can intersect at most
\[  
8M(M+1) Q^{1+2\varepsilon}\cdot 2r+2 \ls 16 (M+1)^2Q^{2\varepsilon}  
\]
intervals centred at points $n(\ba,\bb)\alpha$ of length at most $8 M q_{J}(\ba)^{-1}$. We conclude that for a fixed $\ba$, there are at most $16 (M+1)^2 Q^{2\varepsilon}$ values of $\bb$ such that the set $C(\ba,\bb)$ satisfies~\eqref{lem_FD_nontrivial_intersection}.

Combining the two bounds, we deduce that the number of sets $C(\ba,\bb)$ that satisfy \eqref{lem_FD_nontrivial_intersection} is at most  
\[
5Q^{1+4\varepsilon}\cdot 16 (M+1)^2 Q^{2\varepsilon} = 80 (M+1)^2 Q^{1+6\varepsilon}.
\] 
For any such set the upper bound of Lemma~\ref{lem_size_of_cylinder_set}, together with~\eqref{T_m_is_large}, gives  
\begin{equation} \label{ub_nu_C}
\nu(C(\ba,\bb)) 
< Q^{(1+2\varepsilon)
	(1-2\kappa)}.
\end{equation}
Using~\eqref{ub_nu_C} together with the upper bound for the number of cylinder sets satisfying~\eqref{lem_FD_nontrivial_intersection}, we infer the upper bound
\[
\begin{aligned}
\nu\left(B_r(\bx)\right) &< 80 (M+1)^2 Q^{1+6\varepsilon} \cdot  Q^{(1+2\varepsilon)(1-2\kappa)}  \\
&= 80 (M+1)^2 \cdot  Q^{2-2\kappa+4(2-\kappa)\varepsilon}.
\end{aligned}
\]
Using \eqref{r_level}, this results in the upper bound
\[        \nu\left(B_r(\bx)\right) \le 80 (M+1)^2 Q_0
r^{2\kappa-2-4(2-\kap)\varepsilon}.
\]
\end{proof}

\subsection{Some measure estimates}

In this section, our goal is to prove Proposition \ref{prop_close_ratios} below. Recall the definition of the set $\cE$ in \eqref{Edef}. We additionally define the set
\[
\cE_1 =\left\{\ba\in [M]^{J_0} :  (\ba,\bb)\in \cE  \text{ for some }\bb\in [M]_{0}^{J_0} \right\}.
\]
The set $\cE_1$ can be thought of as the projection of the `two-dimensional' set $\cE$ onto the first coordinate. 

As before, for $N\gs 1$ we write $J=J_0 N$ and 
\[
Q = Q_N = \exp(J\sig_m(\kap)).
\]
Given $\ba = (a_1, \ldots, a_{J}),$ we write $\bar{\ba} = (a_{J}, \ldots, a_1)$ for the mirror image of the word $\ba$. We are interested in sets $A\subseteq \cE^N$ such that \begin{equation} \label{full}
\begin{aligned}
&(\ba, \bb_0) \in A \text{ for some } \bb_0\in \cV_{J}(\ba) \\
&\Rightarrow \quad (\ba, \bb) \in A \text{ for all }  
\bb \in \cV_{J}(\ba).
\end{aligned}
\end{equation}
For any such set $A\subseteq \cE^N,$ we define 
\begin{equation} 
\label{mirror_set} 
\bar{A} = \left\{ (\bar{\ba},\bb): 
\begin{aligned}
&\bb\in \cV_{J}(\bar{\ba}), \\
&(\ba,\bb')\in A \text{ for some } \bb' \in \cV_{J}(\ba)
\end{aligned}
\right\}. 
\end{equation} 

\begin{prop} \label{prop_almost_invariance}
Let $N\gs 1$ and $A\subset\cE^N$ be a set satisfying \eqref{full}. Then 
\[  
\lam^N(\bar{A})
Q^{-(1+2\varepsilon)c/(R+c)}\leq \lam^N(A)
\leq \lam^N(\bar{A})
Q^{(1+2\varepsilon)
c/(R+c)}, 
\]
where the constant $c$ is as in Lemma~\ref{vq}.
\end{prop}

\begin{proof} It is enough to prove the proposition for sets of the form  
\[
A=\{(\ba,\bb): \bb \in \cV_{J}(\ba)\}
\]
for some $\ba=(a_1,\ldots,a_{J})$; the result will then follow  by taking unions over words $\ba \in [M]^{J}$.
Note from the definition of $\cE_1$ that we have $\ba \in \cE_1^N$ if and only if $\bar{\ba}\in \cE_1^N.$ The result is trivially true when $\ba\notin \cE_1^N,$ as then $\nu(A) = \nu(\bar{A})=0.$ 

If $\ba\in \cE_1^N$ then, for $A$ as above, by the definition of $\lam$ we have 
\[
\lam^N(A) = \frac{c_{\lambda}^N |\cV_{J}(\ba)|}{T_m(\kappa)^{j_0 N}}\prod_{i=0}^{j_0 N-1}q_m\left(a_{im+1},\dots, a_{(i+1)m}\right)^{1-2\kappa}.        
\]
In light of the symmetry \eqref{q_k_inverse}, we similarly have
\[
\lam^N(\bar{A}) = \frac{c_{\lambda}^N |\cV_{J}(\bar{\ba})|}{T_m(\kappa)^{j_0 N}}\prod_{i=0}^{j_0 N-1}q_m\left(a_{im+1},\dots, a_{(i+1)m}\right)^{1-2\kappa}.        
\]
By Lemma~\ref{vq}, we find that for every $\ba\in\cE_1^N$,
\[ |\cV_{J}(\ba)|\leq q_{J}(\ba) \ls |\cV_{J}(\ba)|^{1+\frac{c}{R}} \quad \text{ and } \quad |\cV_{J}(\bar{\ba})| \leq q_{J}(\bar{\ba}) \leq |\cV_{J}(\bar{\ba})|^{1+\frac{c}{R}}. \] 
By  \eqref{q_k_inverse} we have $q_{J}(\ba) = q_{J}(\bar{\ba})$ and therefore,  employing~\eqref{Qest}, we get
\[
\begin{aligned}
|\cV_{J}(\ba)| &\gs |\cV_{J}(\bar{\ba})|^{1/(1+\frac{c}{R})} = |\cV_J(\bar{\ba})| \cdot |\cV_{J}(\bar{\ba})|^{-c/(c+R)}\\ &\gs |\cV_{J}(\bar{\ba})| Q^{-(1+2\varepsilon)c/(R+c)}.\end{aligned} 
\]
This gives
\[ 
\lam^N(A) \gs 
\lam^N(\bar{A}) 
Q^{-(1+2\varepsilon)
c/(R+c)},
\]
and the other inequality follows by interchanging the roles of the sets $A$ and $\bar{A}.$ 
\end{proof}

\begin{prop} \label{prop_close_ratios}
Let $N\gs 1$  and $0<T<1$. Define
\[
S = \left\{ (\ba,\bb,\ba^*,\bb^*)\in \cE^N \times \cE^N :  \left|\frac{q_{J-1}}{q_J}-\frac{q_{J-1}^*}{q_J^*}\right|<T \right\}.
\]
Then
\[
\lam^{2N}(S) \ll 
T^{2\kappa-3-4(2-\kap)\varepsilon}
Q^{(1+2\varepsilon)2c/(R+c)},
\]
where the constant $c$ is as in Lemma~\ref{vq}. 
\end{prop}

\begin{proof}
By~\eqref{cf_inverse}, we have
\[  
\left|\frac{q_{J-1}}{q_J}-\frac{q_{J-1}^*}{q_J^*}\right|
= \left|\frac{p_J(\bar{\ba})}{q_J(\bar{\ba})}-\frac{p_J(\bar{\ba^*})}{q_J(\bar{\ba^*})}\right|.
\]
Observe that, upon mapping each $(\ba, \bb, \ba^*, \bb^*) \in S$ to $(\ba \ba^*, \bb \bb^*)$, the set $S$ satisfies the property described in \eqref{full}. The set $\bar{S}$ as in \eqref{mirror_set} is therefore well defined, and
\[
\bar{S} \subseteq \left\{ (\ba,\bb,\ba^*,\bb^*) \in \cE^N\times \cE^N :  \left|\frac{p_J}{q_J}-\frac{p_J^*}{q_J^*}\right| < T \right\},
\]
where $p_J^* = p_J(\ba^*)$ and $q_J^* = q_J(\ba^*)$.

For any $(\ba,\bb)\in \cE^N,$ the set 
\[ 
\bar{S}(\ba,\bb) := \left\{ (\ba^*, \bb^*) \in \cE^{N} : \left| \frac{p_J}{q_J} - \frac{p_J^*}{q_J^*}  \right| < T  \right\} 
\]
lifts to a vertical strip with base length $2T$ given by
\[
\tilde S(\ba, \bb)
= \bigcup_{(\ba^*,\bb^*) \in \bar{S}(\ba,\bb)} 
C(\ba^*, \bb^*).
\]
If we cover it by $1/(2T)$ squares of side $2T$, then each such square is contained in a ball of radius $O(T)$. Therefore
\[  
\lam^N
(\bar{S}(\ba,\bb))
= \nu(\tilde{S}(\ba,\bb))
\ll \frac{1}{2T}\cdot T^{2\kappa - 2 -4(2-\kap)\varepsilon} \ll
T^{2\kappa - 3 - 4(2-\kap)\varepsilon 
} 
\]
and  
\[ 
\lam^{2N}(\bar{S}) \ll \int \lam^N(\bar{S}(\ba, \bb))\d \lam^N \ll T^{2\kappa - 3 - 4(2-\kap)\varepsilon}. 
\]  
The conclusion now follows from Proposition \ref{prop_almost_invariance}.
\end{proof}

\section{Polynomial Fourier decay} \label{section_FD}

In this section, we will upper bound the decay of the Fourier transform of the measure constructed in Section~\ref{subsection_ConstructionOfMeasures}, thereby establishing
Theorem~\ref{MainThm}~\eqref{MainThm_point_decay}. As Theorem~\ref{MainThm}~\eqref{MainThm_point_Frostman_dimension} follows from Lemma~\ref{lem_FD} point~\eqref{point_dimFrostman}, this will complete the proof of Theorem \ref{MainThm}.

\subsection{A comparison lemma} Given some basic data, the following lemma enables us to compare an integral with respect to a general measure to an integral with respect to Lebesgue measure, in two dimensions. This is analogous to its one-dimensional predecessor \cite[Lemma 5.3]{QR2003}.

\begin{lemma} \label{ComparisonLemma}
Let $M_1, M_2 \in (0,\infty)$, let $F(t,u)$ be a $C^1$ function on $B = [0,1] \times [-1,1]$ such that
\[
\| F \|_\infty \le 1,
\qquad 
\| \partial_t F\|_\infty \le M_1, 
\qquad
\| \partial_u F \|_\infty \le M_2,
\]
and put
\[
m_2 = \int_B |F(t,u)|^2 \d t \d u.
\]
Let $\nu$ be a probability measure on $B$, let $r > 0$, and denote by $\Lam(r)$ the maximum of $\nu(\cR)$ as $\cR$ ranges over axis-parallel filled rectangles of respective side lengths $r/M_1$ and $r/M_2$ in $B$. Then
\[
\int_B |F(t,u)| \d \nu(t,u) \le 
3r + \Lam(r) \left(2 + \frac{M_1}r + \frac{2M_2}r + \frac{m_2 M_1 M_2}{r^4} \right).
\]
\end{lemma}

\begin{proof} Consider $\lfloor M_1/r \rfloor \cdot \lfloor 2 M_2/r \rfloor$ axis-parallel filled rectangles of respective side lengths $r/M_1$ and $r/M_2$, with non-overlapping interiors, whose union covers
\[
\left[0, 1 - \: \frac{r}{M_1}
\right] \times
\left[-1, 1 - \: \frac{r}{M_2}
\right].
\]
Let $N$ be the number of these filled rectangles $\cR$ on which
\[
\max \{ |F(t,u)|: (t,u) \in \cR \} \ge 3r.
\]
Within any such filled rectangle $\cR$, there exists a pair $(t_0,u_0)$ such that $|F(t_0,u_0)| \ge 3r$. Thus, for any $(t,u) \in \cR$, we have
\begin{align*}
|F(t,u)| 
&\ge |F(t_0,u_0)| - |F(t,u)-F(t,u_0)| - |F(t,u_0) - F(t_0,u_0)| \\
&\ge 3r - r - r = r,
\end{align*}
by the mean value theorem.
Therefore
\[ 
m_2 \ge N \frac{r^2}{M_1 M_2} r^2 = N \frac{r^4}{M_1M_2}.
\]
At most $\lceil M_1/r \rceil + \lceil 2M_2 / r \rceil$ filled rectangles, further to our original ones, are needed to cover $B$, so
\begin{align*}
\int_B |F(t,u)| \d \nu (t,u) &\le 3r + (N + \lceil M_1/r \rceil + \lceil 2M_2 / r \rceil) \Lam(r) 
\\ &\le 3r + \Lam(r) \left(2 + \frac{M_1}r + \frac{2M_2}r + \frac{m_2 M_1 M_2}{r^4} \right).
\end{align*}
\end{proof}

\subsection{The Fourier transform of 
our measure} Let $k_1, k_2 \in \mathbb{Z}.$ The Fourier transform of $\nu$ evaluated at the point $(k_1, k_2)$ is 
\[ 
\widehat \nu(k_1,k_2) = \int_{\bR^2} e(-k_1 \alp - k_2 \gam) \d \nu(\alp,\gam). \]
Since we cannot derive an explicit formula for the Fourier transform, our aim is to provide an estimate by writing $\widehat{\nu}(k_1,k_2)$ as an integral of the form $\int_B F(t,u)\d \nu(t,u)$ for some appropriate function $F(t,u)$ and then applying Lemma \ref{ComparisonLemma}.

In what follows, we use the parameters
\begin{equation} \label{def_delta_phi}
%\delta=1/5 \text{ and } \phi=\delta/10. 
\delta>0 \text{ and } \phi=\delta/10.
\end{equation}
We also assume that the parameter $\varepsilon>0$ used in the construction of the measure $\nu$, from the beginning of Section~\ref{subsection_ConstructionOfMeasures}, is small enough to satisfy $\varepsilon<\phi/100$.

\begin{remark}
For the proof, we may assume that $\delta$ is large, see Remark~\ref{rem_justification_concrete_eta}. For $\delta$ large enough, some of the cases considered below become simpler, for instance the cases described by~\eqref{third_scenario_description} or~\eqref{Case_2i_subcase_one}. For the sake of robustness, we also allow $\del$ to be small. However, the busy reader may wish to assume that $\del$ and hence $\phi$ is sufficiently large.
For $\phi\geq 2$ (and $Q$ large enough, which we can assume without loss of generality), the inequality ~\eqref{third_scenario_description} implies that
$q_J=q_J^*$ and 
$q_{J-1}=q_{J-1}^*$, and~\eqref{Case_2i_subcase_one} becomes $q_{J-1}=q_{J-1}^*$. The two simplified cases then reduce to the diagonal setting considered at the start of Subsection~\ref{ss_diagonal_contribution}.
\end{remark}

With the pair $(k_1, k_2)$ fixed, let $K \ge 0$ be an even integer such that
\[
Q_0^{KJ_0(2+\delta)}
\leq\max
\left(|k_1|, |k_2| 
Q_0^{KJ_0}\right)\leq Q_0^{(K+2)J_0(2+\delta)}.
\]
Also set
\[
J := KJ_0 \qquad \text{ and } \qquad
Q := Q_0^{J}.
\]
It then follows that
\begin{equation} \label{kk_size}
\max\left(|k_1|,|k_2| Q\right)\asymp Q^{2+\delta}.
\end{equation}
Writing
\[
\alp = [a_1,a_2,\ldots],
\qquad 
\gam = \sum_{k=0}^\infty b_{k+1} D_k(\alp),
\]
as well as
\[
\ba = (a_1,\ldots,a_J), \qquad
\bb = (b_1,\ldots,b_J),
\]
we now have
\begin{equation} \label{factor}
\d \nu(\alp,\gam) = \d \lam^K(\ba,\bb) \d \nu(t,u),
\end{equation}
where
\[
t = [a_{J+1},a_{J+2},\ldots],
\qquad
u = \sum_{k=0}^\infty b_{J+k+1} D_k(t).
\]

\begin{lemma} We have
\[
-t \le u \le 1-t.
\]
\end{lemma}

\begin{proof} Observe that $u$ is given by an Ostrowski expansion in $t$. We use the fact that $D_0, D_1, D_2, \ldots$ alternate in sign, together with the identity \eqref{EndpointIdentity}. Further, as $R \mid J$, we have $b_{J+1} \le a_{J+1} - 1 = a_1(t) - 1$. Using this information, we compute that
\[
u \ge a_2(t) D_1(t) + a_4(t) D_3(t) + \dots
= a_2(t) D_1(t) - D_2(t) = - t
\]
and
\begin{align*}
u &\le (a_1(t) - 1) D_0(t) + a_3(t) D_2(t) + a_5(t) D_3(t) + \dots
\\ &= (a_1(t) - 1) D_0(t) - D_1(t) = 1 - t.
\end{align*}
\end{proof}

Let us now fix
\begin{equation} \label{def_B}
\begin{aligned}
P = \left\{(t,u)\in [0,1]\times\mathbb{R}: -t\leq u\leq 1-t \right\}\subset [0,1]\times[-1,1].
\end{aligned}
\end{equation}

\begin{cor} The pushforward measure $\nu$ is supported on $P$.
\end{cor} 

\begin{remark} \label{period_nu}
Clearly, translations by $\bZ^2$ of $P$, defined by~\eqref{def_B} above, tile the plane $\bR^2$, hence we can extend the measure $\nu$ to $\bR^2$ by periodicity. We could also consider this measure as being defined on $[0,1]^2$. Indeed, the parallelogram $P$ consists of two triangles:
\[
\begin{aligned}
P&=T_1\cup T_2,\\
T_1&:=\left\{(t,u)\in [0,1]\times [0,1] : 0 \leq u\leq 1-t \right\},\\
T_2&:=\left\{(t,u)\in [0,1]\times [-1,0] : -t \leq u\leq 0 \right\},
\end{aligned}
\]
and $[0,1]^2$ is the union of $T_1$ and the translation of $T_2$ by the vector $(0,1)$.
\end{remark}

We now apply the analysis of \S \ref{Ostrowski_section} to describe all quantities associated with pairs of the form $(\alp,\gam)$ as functions of  $(\ba,\bb)$ and $(t,u)$ and write $\widehat\nu(k_1,k_2)$ in the desired form. It will be convenient to write
\[
g = g(\ba,\bb,t)= \sum_{k=0}^{J-1} b_{k+1} D_k(\alp).
\]
In order to view $g$ as a function of $\ba,\bb,t$, we observe that
\begin{align}
&\gamma = g  - D_{J-1}(\alp) u, \label{eq2}\\
&D_k(\alp) = q_k \frac{p_J + tp_{J-1}}{q_J + tq_{J-1}} - p_k \qquad (k \ge 0), \label{eq3} \\
&D_{J-1}(\alp)  = \frac{(-1)^{J+1}}{q_J + tq_{J-1}}. \label{eq4}
\end{align}
Indeed, using \eqref{DkRelation}, we get 
\[
\gamma = \sum_{k=0}^{J-1}b_{k+1}D_{k}(\alp) + \sum_{k=0}^{\infty}b_{J+k+1}D_{J+k}(\alp) = g - D_{J-1}(\alp) u,
\] 
which is \eqref{eq2}. Equation \eqref{eq3} follows from \eqref{alphaoft} and the definition of $D_k$, while \eqref{eq4} follows from \eqref{eq3} and \eqref{pkqk}. 

By \eqref{factor} and the fact that $J$ is even, we now have 
\begin{align*}
\widehat \nu(k_1,k_2)
&= \int_{\bR^2} e(-k_1 \alp - k_2 \gam) \d \nu(\alp,\gam)\\ 
&= \int_{\bR^2} \int_{\cE^K}e \left( -k_1
\frac{p_J + t p_{J-1}}{q_J + tq_{J-1}} \: - k_2 \left(g(\ba,\bb,t) + \frac{u}{q_J + tq_{J-1}}\right) \right) \\
& \qquad \d\lam^K(\ba,\bb) \d\nu(t,u)  \\
&= \int_B F(t,u)  \d \nu(t,u),
\end{align*}
where
\begin{align} 
\notag F(t,u)
&= \int_{\cE^K} e \left( -k_1
\frac{p_J + t p_{J-1}}{q_J + tq_{J-1}} \: - k_2
\left(g(\ba,\bb,t) + \frac{u}{q_J + tq_{J-1}}\right) \right) \\
\label{F_definition}
& \qquad \d \lam^K(\ba,\bb).
\end{align}
By Lemma \ref{ComparisonLemma}, to estimate $\widehat \nu(k_1,k_2)$ it now remains to bound the partial derivatives and second moment of our function $F$, as well as the measures of certain rectangles.

\subsection{The partial derivatives} \label{partial_d}

For $k \ge 0$, we have 
\[
\partial_t\alpha = \partial_t \left( \frac{p_J + t p_{J-1}}{q_J + tq_{J-1}} \right) = \frac{1}{(q_J+tq_{J-1})^2}
\]
and
\[
\partial_t D_k(\alp) = q_k\cdot \partial_t \alp = \frac{q_k}{(q_J+tq_{J-1})^2}.
\]
Hence
\begin{equation} \label{gt}
\partial_t g = \sum_{k=0}^{J-1}b_{k+1}\partial_t D_k(\alp) = (q_J+tq_{J-1})^{-2}\sum_{k=0}^{J-1}b_{k+1}q_k   
\end{equation}
and
\[ 
\partial_t \left( \frac{u}{q_J + tq_{J-1}}\right) = \frac{-u q_{J-1}}{(q_J+tq_{J-1})^2 }.
\]
The partial derivative of $F$ with respect to $t$ is therefore
\begin{align*}
&\partial_t F(t,u) \\
&=  \int_{\cE^K} \frac{-2 \pi i}{(q_J+tq_{J-1})^2}\left(k_1 + k_2 \left(-uq_{J-1} + \sum_{k=0}^{J-1}b_{k+1}q_k \right) \right) \\
&  \qquad \cdot e \left( -k_1 \frac{p_J + t p_{J-1}}{q_J + tq_{J-1}} - k_2 \left(g(\ba,\bb,t) + \frac{u}{q_J + tq_{J-1}}\right) \right) \d \lam^K(\ba,\bb).
\end{align*}
By uniqueness of the Ostrowski expansion, if $(\ba,\bb) \in \cE^K$ then
\[
\sum_{k=0}^{J-1}b_{k+1}q_k <  q_J,
\]
see \cite[Chapter II, \S 4]{RS1992}, and so
\[ 
\left| -uq_{J-1} + 
\sum_{k=0}^{J-1}b_{k+1}q_k 
\right| 
< q_J + q_{J-1} \qquad (-1 \le u \le 1).
\]
Therefore
\[
|\partial_t F(t,u)| \le 2 \pi \left( \frac{|k_1|}{Q^{2(1-2\varepsilon)}} + \frac{2|k_2|}{Q^{1-2\varepsilon}} \right).
\]
The partial derivative of $F(t,u)$ with respect to $u$ is
\begin{align*}
&\partial_u F(t,u) \\
&=  \int_{\cE^K}\frac{-(2 \pi i) k_2}{q_J+tq_{J-1}} e\left( -k_1
\frac{p_J + t p_{J-1}}{q_J + tq_{J-1}} - k_2
\left(g(\ba,\bb,t) + \frac{u}{q_J + tq_{J-1}}\right) \right) \\
&\qquad \d \lam^K(\ba,\bb),
\end{align*}
whence
\[
|\partial_u F(t,u)| \le \frac{2\pi |k_2|}{ Q^{1-2\varepsilon}}.
\]

\subsection{The second moment}
\label{m2calc}

With $B = [0,1] \times [-1, 1]$, the second moment of the function $F(t,u)$ defined in \eqref{F_definition} is 
\begin{align*}
m_2 &= \int_B |F(t,u)|^2 \mathrm{d}t \mathrm{d}u \\ 
&= \int_{\cE^{2K}} \int_0^1 \int_{-1}^1 e\Big(k_1\Big( \frac{p_J+tp_{J-1}}{q_J+tq_{J-1}} \: - \: \frac{p_J^* +tp_{J-1}^*}{q_J^* +tq_{J-1}^*}\Big)\Big)  \\ &
\quad e\big(k_2(g(\ba,\bb,t) - g(\ba^*,\bb^*,t))\big) \\
& \quad e\Big(k_2u\Big( \frac{1}{q_J+tq_{J-1}} \: - \: \frac{1}{q_J^*+tq_{J-1}^*}\Big)\Big)
\d u \d t \d\lambda^K(\ba,\bb) \d\lambda^K(\ba^*, \bb^*).
\end{align*}
Here and in what follows, for $(\ba,\bb), (\ba^*,\bb^*) \in \cE^K$, we write for convenience  $q_J=q_J(\ba)$ and $q_J^*=q_J(\ba^*),$ and similarly for $q_{J-1}, q_{J-1}^*.$ 

Continuing with the estimation of $m_2,$ let us set
\begin{equation}\label{f1f2} f_1(t) = \frac{p_J + t p_{J-1}}{q_J + tq_{J-1}} \: - \: \frac{p_J^* + tp_{J-1}^*}{q_J^* + tq_{J-1}^*},
\qquad
f_2(t) = g(\ba,\bb, t) - g(\ba^*, \bb^*, t)
\end{equation}
and
\begin{equation} \label{def_z}
z(t) = z(t;k_2) = k_2 \left( \frac1{q_J + tq_{J-1}} \: - \: \frac1{q_J^* + t q_{J-1}^*} \right),
\end{equation}
and write
\begin{equation} \label{f_def}
f(t) = f(t; k_1, k_2) = k_1 f_1(t) + k_2 f_2(t).
\end{equation}
Now
\begin{align} \label{m2} 
m_2 &= \int_{\cE^{2K}} \fJ \d \lam^{2K},
\end{align}
where
\[
\fJ = \fJ(\ba, \bb, \ba^*, \bb^*) = \int_0^1 e(f(t)) \int_{-1}^1 e(u z(t)) \d u \d t.
\]
Observe that if $z \in \bR$ then
\[
\int_{-1}^1 e(zu) \d u = 2 \sinc(2 \pi z).
\]
This is immediate if $z = 0$, and otherwise
\[
\int_{-1}^1 e(zu) \d u  
= \frac{e(z) - e(-z)}{2 \pi i z} 
= \frac{\sin (2 \pi z)}{\pi z} 
= 2 \sinc(2\pi z).
\]
Therefore
\begin{equation} \label{fJ_first_presentation}
\fJ = \fJ(\ba, \bb, \ba^*, \bb^*) = 2 \int_0^1 e(f(t)) \sinc(2 \pi z(t)) \d t.
\end{equation}
We now estimate $m_2$ by partitioning $\cE^K\times \cE^K$ into suitably chosen subsets and bounding the corresponding contribution of $\fJ$. 

\subsubsection{The diagonal contribution} \label{ss_diagonal_contribution}
Write 
\[
\cD = \{ (\ba, \bb, \ba^*, \bb^*) \in \cE^{2K}: \ba = \ba^* \}
\]
for the set of quadruples $(\ba,\bb,\ba^*,\bb^*)\in \cE^{2K}$ such that $\ba = \ba^*$. Since $|\fJ| \le 2,$ the contribution of $\cD$ to the integral in \eqref{m2} is
\begin{equation} \label{diagonal}
\int_{\cD}\fJ \d \lam^{2K} \ll   (\lambda^K \times \lambda^K)(\cD).
\end{equation}   
The right-hand side of \eqref{diagonal} is equal to 
\begin{align*}
(\lambda^K \times \lambda^K)(\cD) & = \sum_{\ba \in \cE_1^{K}}\mathop{\sum\sum}_{\bb, \bb^* \in \cV_{KJ_0}(\ba)}(\lambda^K \times \lambda^K)(\ba,\bb,\ba,\bb^*)\\
& = \sum_{\ba \in \cE_1^{K}} \Bigg(\sum_{\bb\in \cV_{KJ_0}(\ba)}\!\!\lambda^K(\ba,\bb)\Bigg)^2 .
\end{align*}
By Lemmas \ref{lem_size_of_cylinder_set} and \ref{vq}, the innermost sum   satisfies
\begin{align*}
\sum_{\bb\in \cV_{KJ_0}(\ba)}\!\!\lambda^K(\ba,\bb) 
&\ls 
| \cV_{KJ_0}(\ba)| \frac{2^K}{T_m(\kappa)^{j_0K}}Q_0^{(1+2\varepsilon)(1-2\kappa)K}\\
&\le \frac{2^{K}Q_0^{(1-2\kap)(1+2\eps)K+(1+2\eps)K}}{T_m(\kap)^{j_0K}}.
\end{align*}
This gives
\begin{align*}
(\lambda^K \times \lambda^K)(\cD) &\ls \frac{2^{K}Q_0^{(1-2\kap)(1+2\eps)K+(1+2\eps)K}}{T_m(\kap)^{j_0K}} \sum_{\ba \in \cE_1^{K}}  \sum_{\bb\in \cV_{KJ_0}(\ba)}\!\!\lambda^K(\ba,\bb) \\
& = \frac{2^{K}Q_0^{(1-2\kap)(1+2\eps)K+(1+2\eps)K}}{T_m(\kap)^{j_0K}},
\end{align*}   
where in the last step we have used the fact that $\lambda^K$ is a probability measure. The assumption \eqref{T_m_is_large} now yields 
\[   
(\lambda^K \times \lambda^K)(\cD)  \ls  Q_0^{(1-2\kap)(1+2\eps)K+(1+2\eps)K}  
\]
and, recalling that $Q = Q_0^K$, we finally have
\[  
\int_{\cD}\fJ \d \lam^{2K}  
\ll  Q^{2-2\kappa 
+ 4(1-\kappa)\varepsilon}.  
\]

\subsubsection{The off-diagonal contribution} We now turn to the contribution of the complement of $\cD$ to the integral in \eqref{m2}. Bounding the $\lam^{2K}$ measure of the domain of integration alone will not be enough, so we aim to establish alternative upper bounds for the integrand $\fJ.$ \par

Write 
\begin{equation} \label{sha_def}
h(t) = 2 \sinc(2 \pi z(t)),
\quad \fA(y) = \int_0^y e(f(t)) \d t \qquad (0\ls t,y \ls 1).
\end{equation}
By partial integration in \eqref{fJ_first_presentation}, we have
\[
\fJ = \int_0^1 e(f(t)) h(t) \d t = \fA(1) h(1) - \int_0^1 \fA(y) h'(y) 
\d y.
\]
Since for the function $ z(t)$ defined in \eqref{def_z} we have $z'(t) \ll |k_2| Q^{2\eps - 1}$, and $\sinc$ has bounded derivative, we have
\[
h'(t) \ll |k_2| Q^{2\eps - 1} 
\]
and consequently
\begin{equation} \label{ub_J_second}
\fJ \ll (1 + |k_2| Q^{2\eps - 1}) \max \{ |\fA(y)|: 0 \le y \le 1 \}.
\end{equation}

We bound the off-diagonal contribution to $m_2$ by partitioning $\cE^{2K}\setminus \cD$ into appropriately chosen subsets, each of them corresponding to one of the cases that follow.

\bigskip

%{\color{red} This case essentially disappears when $\delta$ is large enough}

\textbf{Case~I.} Recalling the parameter $\phi>0$ from~\eqref{def_delta_phi}, we consider the quadruples $(\ba,\bb,\ba^*,\bb^*)$ for which
\begin{equation} \label{third_scenario_description}
|q_J - q_J^*| +
|q_{J-1} - q_{J-1}^*| 
\le 3Q^{1-\phi}.
\end{equation}
Under this assumption, we compute that
\[
\begin{aligned}
\left|\frac{q_{J-1}}{q_J}-\frac{q_{J-1}^*}{q_J^*}\right|
&= \frac{|q_{J-1} q_J^*-q_{J} q_{J-1}^*|}{q_{J} q_J^*} \\
&= \frac{\left|q_{J-1} \left(q_J^*-q_J\right)-q_{J}\left(q_{J-1}^*-q_{J-1}\right)\right|}{q_{J}q_J^*}\\
&\ll \frac{|q_J^* - q_J| 
+ |q_{J-1}^* - q_{J-1}|}{q_J^*} \ll Q^{2\varepsilon-\phi},
\end{aligned} 
\]
where at the last step we used Lemma~\ref{lem_FD}~(a) and~\eqref{third_scenario_description}. Thus, by Proposition~\ref{prop_close_ratios},
we have
\begin{align*}
&(\lambda^K \times \lambda^K)(\{ (\ba,\bb,\ba^*,\bb^*) :  \text{ \eqref{third_scenario_description} holds}\})
\\ & \qquad \ll  Q^{(2\eps - \phi)
(2 \kap - 3 - 4(2-\kap)\eps)
+ (1 + 2\eps) 2c / (R+c)} ,
\end{align*}
where the constant $c$ is as in Lemma~\ref{vq}. We conclude that the  contribution to $m_2$ is 
\begin{equation}
\label{Contribution1}
O( Q^{(2\eps - \phi)
(2 \kap - 3 - 4(2-\kap)\eps)
+ (1 + 2\eps) 2c / (R+c)} ). 
\end{equation}
It remains to deal with the contribution to $m_2$ of $(\ba,\bb,\ba^*,\bb^*)$ for which \eqref{third_scenario_description} fails, and we do so by considering two subcases according to the size of $|k_2|.$ 

\bigskip 

\textbf{Case~II(i).} Here we examine the contribution to $m_2$ of $(\ba,\bb,\ba^*,\bb^*)$ such that 
\begin{equation}\label{sum_q_big}
|q_J^* - q_J| +
|q_{J-1}^* - q_{J-1}| 
> 3Q^{1- \phi} 
\end{equation}
when
\begin{equation} \label{K_is_big}
|k_2| \gg Q^{1 + 3\phi}.
\end{equation}
By the definitions of $\fJ$ and $h(t)$ in \eqref{fJ_first_presentation}, \eqref{sha_def}, we have
\[ \fJ \ll  \int_0^1 |h(t)|\d t \ll  \int_0^1 \min
\{ 1, |z(t)|^{-1} \} \d t, \]
where by \eqref{def_z}
\[
z(t) = k_2
\frac{t(q_{J-1}^* - q_{J-1})
+ q_J^* - q_J}
{(q_J + tq_{J-1})
(q_J^* + tq_{J-1}^*)}.
\]

First, suppose
\begin{equation} \label{FirstSuppose}
|q_{J-1}^* - q_{J-1}| > Q^{1-\phi}.
\end{equation}
Write $S$ for the set of $t \in [0,1]$ such that
\[
\left| t + 
\frac{q_J^* - q_J}
{q_{J-1}^* - q_{J-1}}
\right| \le Q^{- \phi}.
\]
Then
\[
\fJ \ll \fJ_1 + \fJ_2,
\]
where
\[
\fJ_1 = \int_S 
\min
\{ 1, |z(t)|^{-1} \} \d t,
\qquad
\fJ_2 = \int_{[0,1]\setminus S}
\min
\{ 1, |z(t)|^{-1} \} \d t.
\]
For $\fJ_1,$ we bound the integrand by $1$ and deduce that $\fJ_1 \ll Q^{-\phi}.$ Regarding $\fJ_2,$ the inequalities \eqref{denominator_bounds},
\eqref{K_is_big} and \eqref{FirstSuppose} imply that for any $t\in [0,1] \setminus S$ we have 
\[
|z(t)| \gg \: Q^{3\phi - 1 - 4\eps}
|q_{J-1}^* - q_{J-1}| \cdot
\left| t + 
\frac{q_J^* - q_J}
{q_{J-1}^* - q_{J-1}} 
\right| \gg Q^{\phi - 4\varepsilon}. 
\]
Thus, we compute that
$ \fJ_2 \ll Q^{4 \eps - \phi},$
whence $\fJ \ll   Q^{4 \eps - \phi}.$

%{\color{red} This case essentially disappears when $\delta$ is large enough}

Now suppose instead that
\begin{equation} \label{Case_2i_subcase_one}
|q_{J-1}^* - q_{J-1}| \le Q^{1 - \phi}.
\end{equation}
Then by  \eqref{sum_q_big}
\[
|q_J^* - q_J| > 2Q^{1 - \phi},
\]
and consequently for any $0 \le t \le 1,$
\[ \begin{split}
|z(t)| &\gg  Q^{3\phi - 1 - 4 \eps} |(q_{J-1}-q_{J-1}^*)t +(q_J-q_J^*)|
> Q^{2 \phi - 4\eps}.
\end{split}
\]
Therefore 
$\fJ \ll Q^{ 4\eps - 2 \phi}$.

We have demonstrated that $\fJ \ll Q^{4\eps - \phi}$ in Case II(i). This means that the corresponding contribution to $m_2$ is 
\begin{equation} 
\label{Contribution21}
O(Q^{4 \eps - \phi}).
\end{equation}

\bigskip

\textbf{Case II(ii).} It remains to deal with the contribution of points for which we have \eqref{sum_q_big} when
\begin{equation} \label{scenario_2}
|k_2| \ll 
Q^{1 + 3 \phi}.
\end{equation}
By \eqref{kk_size} and \eqref{scenario_2}, we have 
\begin{equation}
\label{k1large}
|k_1| \asymp Q^{2 + \del},
\end{equation}
and by~\eqref{ub_J_second} we have
\[
\fJ \ll Q^{3 \phi+2\varepsilon} 
\max \{ |\fA(y)|: 0 \le y \le 1 \}.
\]
For the functions $f_1$ and $f_2$ defined in \eqref{f1f2}, we compute that
\[
f_1'(t) = (q_J+tq_{J-1})^{-2} - (q_J^*+tq^*_{J-1})^{-2}
\]
and, using \eqref{gt},
\begin{align*}
f_2'(t) &=
(q_J+tq_{J-1})^{-2}\sum_{k=0}^{J-1}b_{k+1}q_k
- (q_J^*+tq^*_{J-1})^{-2}\sum_{k=0}^{J-1}b^*_{k+1}q_k^*.
\end{align*}
Thus, for $f=k_1f_1 + k_2f_2$ introduced in \eqref{f_def}, we have
\begin{align*}
f'(t) &=
k_1(
(q_J+tq_{J-1})^{-2} - (q_J^*+tq^*_{J-1})^{-2}) \\
& \quad +
k_2\left( (q_J+tq_{J-1})^{-2}\sum_{k=0}^{J-1}b_{k+1}q_k
- (q_J^*+tq^*_{J-1})^{-2}\sum_{k=0}^{J-1}b^*_{k+1}q_k^* \right)
\\
&=
(q_J+tq_{J-1})^{-2}
\left(k_1+k_2\sum_{k=0}^{J-1}b_{k+1}q_k\right) \\
& \qquad - (q_J^*+tq^*_{J-1})^{-2}
\left(k_1+k_2\sum_{k=0}^{J-1}b_{k+1}^*q_k^*\right).
\end{align*}
Setting
\begin{equation} \label{ss_def}
s = k_1 + k_2 \sum_{k=0}^{J-1} b_{k+1} q_k,
\qquad s^* = k_1 + k_2 \sum_{k=0}^{J-1} b_{k+1}^* q_k^*,
\end{equation}
the derivative $f'(t)$ can be written as
\begin{equation} \label{fprime1}
f'(t) =
\frac{(q_J^* + tq_{J-1}^*)^2 s - (q_J + tq_{J-1})^2 s^*}{(q_J + tq_{J-1})^2 (q_J^* + tq_{J-1}^*)^2}.
\end{equation}
We partition the set 
$\cC :=
\{ (\ba,\bb,\ba^*,\bb^*)\in\cE^{2K}: \eqref{sum_q_big} \text{ holds} \}$ into three sets $\cC_0$, $\cC_1$ and $\cC_2$ that are defined as follows: 
\begin{align*}
\cC_0 &= \{ (\ba,\bb,\ba^*,\bb^*)\in \cC : ss^* =0 \}, \\
\cC_1 &= \{  (\ba,\bb,\ba^*,\bb^*)\in \cC : ss^* <0  \},   \\
\cC_2 &= \{  (\ba,\bb,\ba^*,\bb^*)\in \cC : ss^* >0  \}.
\end{align*}

For the analysis of the sets $\cC_0, \cC_1, \cC_2$, we need to estimate the 
$\lam^K$-measure of the set where $|s|$ assumes values up to a certain threshold. In the following lemma, we essentially consider $s$ as a random variable, where $(\ba,\bb) \in \cE^{K}$ is distributed according to  $ \lam^{K}$.  
\begin{lemma} \label{SmallProb} Assume that $k_2 \ne 0$. Then, for $T \ge 0$, we have
\[
\lam^K(|s| \le T) \le \left(\frac{2T}{|k_2|} + 1 \right) Q^{2 \eps - 1},
\]
where $Q=\exp\left(J\sigma_m\right)$.
\end{lemma}

\begin{proof} First, note that
\[
\lam^K\left(|s| \le T\right)=\bE\left(\chi_{|s| \le T}\right),
\]
where $\chi_{|s| \le T}$ is a characteristic function of the subset of $\cE^K$ defined by the condition $|s| \le T$. Then, by using the tower property of conditional expectation, we have
\[
\bE\left(\chi_{|s| \le T}\right)=\bE\left(\bE\left(\chi_{|s| \le T}\mid \ba\right)\right).
\]
We are going to upper bound the random variable $\bE\left(\chi_{|s| \le T}\mid \ba\right)$, and then use the fact that an upper bound for a random variable is also an upper bound for its expectation.

To find a suitable upper bound for $\bE\left(\chi_{|s| \le T}\mid \ba\right)$, we observe first that
\[
s = k_1 + k_2 N, \qquad N = \sum_{k=0}^{J-1} b_{k+1} q_k.
\]
Writing $\ba = (\ba^{(1)}, \ldots,\ba^{(J/m)})$, the conditional probability that $N$ takes a specific value is at most $\prod_{i=1}^{J/m} |\cV_m(\ba^{(i)})|^{-1}$.  As $s \equiv k_1 \mmod |k_2|$, the conditional probability that $|s| \le T$ is at most
\[
\bE\left(\chi_{|s| \le T}\mid \ba\right) \le \left( \frac{2T}{|k_2|} + 1 \right) \prod_{i=1}^{J/m} |\cV_m(\ba^{(i)})|^{-1}.
\]
By Lemma~\ref{vq} and the inequalities \eqref{R_is_large}, \eqref{Edef}, we infer that
\[
\begin{aligned}
\bE\left(\chi_{|s| \le T}\mid \ba\right) &\le \left(
\frac{2T}{|k_2|} + 1 \right)\frac{27^{J/R}} {\prod_{i=1}^{J/m}q_m(\ba^{(i)})} \\
&\le \left(
\frac{2T}{|k_2|} + 1 \right)\frac{27^{J/R}}{Q^{1-\varepsilon}}
\le \left(
\frac{2T}{|k_2|} + 1 \right)Q^{2\varepsilon-1}.
\end{aligned}
\]
Taking expectations gives the conclusion of the lemma.
\end{proof}

We may now examine each of the sets $\cC_0, \cC_1, \cC_2$ separately.  By Lemma~\ref{SmallProb}, we have $\lam^{2K}(\cC_0) \le 2Q^{2\eps-1}$, and since $|\fJ| \le 2$, this gives
\begin{equation}
\label{Contribution220}
\int_{\cC_0} \fJ \d \lam^{2K} \ll Q^{2\eps - 1}.
\end{equation}
Observe that if $(\ba,\bb,\ba^*,\bb^*) \in \cC_1$ then, by \eqref{fprime1},
\[ Q^{2+12\eps} |f'(t)| \gg |s| + |s^*| \ge 1. \]
As $f'$ is monotonic, Lemma \ref{lem_ub_nozero} yields, for every $y\in[0,1]$,
\begin{equation} \label{ub_U_C_one}
|\fA(y)(\ba,\bb,\ba^*,\bb^*)| \ll \frac{Q^{2+12\eps}} {|s|+|s^*|}.
\end{equation}
Under the assumption~\eqref{scenario_2}, we have $|k_1|\asymp Q^{2+\delta}$ and thus
\[
|k_1| \gg  |k_2| Q^{1 + 2 \eps}.
\]
Then, by the triangle inequality for \eqref{ss_def}, we have
\begin{equation} \label{lb_s_C_one}
\min\left(|s|, |s^*|\right) \ge |k_1| - |k_2| Q^{1+2\eps} \gg Q^{2+\del}.
\end{equation}
Thus \eqref{ub_J_second}, \eqref{scenario_2}, \eqref{ub_U_C_one} and~\eqref{lb_s_C_one} imply that 
\[
\fJ \ll Q^{3\phi+2\varepsilon-\delta+12\varepsilon}
=Q^{3\phi-\delta+14\varepsilon}
\qquad 
((\ba,\bb,\ba^*,\bb^*) \in \cC_1).
\]
Hence
\begin{equation} 
\label{Contribution221}
\int_{\cC_1} \fJ \d \lam^{2K} \ll Q^{3\phi-\delta+14\varepsilon}.
\end{equation}

\bigskip

Finally, we turn our attention to $\cC_2$. We assume that $s > s^* > 0$. The other cases within $\cC_2$ can be handled by symmetrical reasoning. By \eqref{fprime1}, we have that $f'(t)$ is equal to 
\begin{align*}
&  \frac{((q_J^* + tq_{J-1}^*)v - (q_J + tq_{J-1}) v^*) ((q_J^* + tq_{J-1}^*)v + (q_J + tq_{J-1}) v^*)}{(q_J + tq_{J-1})^2 (q_J^* + tq_{J-1}^*)^2},
\end{align*}
where 
\[
v = \sqrt s, \qquad v^* = \sqrt{s^*}.
\]
Hence
\[
f'(t) = (At + B) G(t),
\]
where
\[
A = q_{J-1}^*v - q_{J-1}v^*, \qquad
B = q_J^* v - q_J v^*
\]
and
\[
G(t) = \frac{(q_J^* + tq_{J-1}^*)v + (q_J + tq_{J-1}) v^*}{(q_J + tq_{J-1})^2 (q_J^* + tq_{J-1}^*)^2}.
\]
Plainly $$|G(t)| \gg Q^{- 3-10\varepsilon} (v+v^*).$$ 
We compute that 
\begin{multline*}
(q_J + tq_{J-1})^3 (q_J^* + tq_{J-1}^*)^3 G'(t)  =  (q_J + tq_{J-1}) (q_J^* + tq_{J-1}^*) (q_{J-1}^*v + q_{J-1} v^*)  \\
\hfill -2(q_{J-1}(q_J^* + tq_{J-1}^*) + q_{J-1}^* (q_J + tq_{J-1})) \cdot \\
\cdot ((q_J^* + tq_{J-1}^*)v + (q_J + tq_{J-1})v^*),
\end{multline*}
and so $$|G'(t)| \ll Q^{ - 3 + 18\eps}(v+v^*).$$ Now for the function $\fA$ defined in \eqref{sha_def},  Lemma \ref{ExpInt} furnishes
\begin{equation} \label{bound_case_one}
\fA(y)^2 \ll \frac{Q^{3+66\eps}} {|A| (v+v^*)} \ll  \frac{Q^{4+68\eps}}{|(q_{J-1}^*)^2 s - q_{J-1}^2 s^*|}   \qquad (0\ls y \ls 1). 
\end{equation}

In order to make use of the bound \eqref{bound_case_one}, we write the region of integration $\cC_2$ as the union of the three sets 
\[ 
\begin{aligned}
&\cS_1 = \left\{
|(q_{J-1}^*)^2 s - q_{J-1}^2 s^*|\geq Q^{4+8\phi} \right\},\\
&\cS_2 = \left\{
\min_{t\in [0,1]}\left|(q_J^* + t q_{J-1}^*)^2 s - (q_J + t q_{J-1})^2 s^*\right|\geq Q^{4+8\phi} \right\}, \\
&\cS_3  
= \left\{
\begin{array}{l}
|(q_{J-1}^*)^2 s - q_{J-1}^2 s^*|<Q^{4+8\phi} \quad \text{ and } \\ \displaystyle \min_{t\in [0,1]}\left|(q_J^* + t q_{J-1}^*)^2 s - (q_J + t q_{J-1})^2 s^*\right|< Q^{4+8\phi}  
\end{array}\right\}.
\end{aligned}
\]
$\bullet$ For $(\ba,\bb,\ba^*,\bb^*)\in\cS_1,$ the bound~\eqref{bound_case_one} gives $\fA(y)\ll Q^{34\eps - 4\phi}$ for every $y\in[0,1]$, and thus, by \eqref{ub_J_second} and \eqref{scenario_2},
\begin{equation}
\label{Contribution2221}
\int_{\cS_1}\fJ \d \lam^{2K} \  \ll Q^{36\eps - \phi}.
\end{equation}
$\bullet$ When $(\ba,\bb,\ba^*,\bb^*)\in\cS_2$, it follows from~\eqref{fprime1} and the definition of $\cS_2$ that  $|f'(t)|\geq Q^{8\phi - 8\varepsilon}$. We also compute, by differentiating~\eqref{fprime1}, that
\[
f''(t) = - \: \frac{2q_{J-1}s}{(q_J + tq_{J-1})^3} +
\frac{2q_{J-1}^*s^*}{(q_J^* + tq_{J-1}^*)^3}.
\]
This shows that $f''(t)=0$ if and only if
\[
q_J^* + tq_{J-1}^*=\left(q_J + tq_{J-1}\right)\left(\frac{q_{J-1}^*s^*}{q_{J-1}s}\right)^{1/3}.
\]

This equation is at most linear in $t$, and we claim that it has at most one solution. Indeed, assume for a contradiction that it has more than one solution. Then $f'' = 0$, so $f'$ is constant. By \eqref{fprime1}, this entails that $f' = 0$, contradicting the definition of $\cS_2$.

Therefore $f'$ is monotonic in at most two stretches. We may therefore apply Lemma~\ref{lem_ub_nozero} with $k=1$ and $\lambda=Q^{8\phi - 8\varepsilon}$, to find
\[
\fA(y)\ll Q^{-8\phi + 8\varepsilon}.
\]
Combining this with \eqref{ub_J_second} and~\eqref{scenario_2} furnishes
\begin{equation}
\label{Contribution2222}
\int_{\cS_2}\fJ \d \lam^{2K} 
\  \ll Q^{10\eps-5\phi}.
\end{equation}
$\bullet$ For $(\ba,\bb,\ba^*,\bb^*)\in\cS_3$, we use the trivial upper bound $|\fJ| \le 2$, and establish an upper bound for $\lambda^{2K}(\cS_3).$ Recall that we have
\eqref{k1large} in the present setting.
Now, applying the triangle inequality to the definitions of $s, s^*$ in \eqref{ss_def}, we have 
\begin{equation} \label{lb_s}
|s|, |s^*| \ge |k_1| - |k_2| Q^{1+2\eps} \gg Q^{2+\del}.
\end{equation}
By definition, for every quadruple in $\cS_3$, we have the following system of linear inequalities:
\begin{equation} \label{system_one}
\begin{cases}
&|({q_{J-1}^*})^2 s - q_{J-1}^2 s^*|<Q^{4+8\phi},\\
&\left|(q_J^* + t q_{J-1}^*)^2 s - (q_J + t q_{J-1})^2 s^*\right| < Q^{4+8\phi}.
\end{cases}
\end{equation}
Since $v=\sqrt{s}$ and $v^*=\sqrt{s^*}$ by definition, we see from \eqref{lb_s} that
\begin{equation} \label{lb_v}
v, v^* \gg Q^{1+\delta/2}.
\end{equation}
Furthermore, we have
\[
\begin{aligned}
(q_{J-1}^*)^2 s - q_{J-1}^2 s^*&=\left(q_{J-1}^* v - q_{J-1} v^*\right)\left(q_{J-1}^* v + q_{J-1} v^*\right)\\
(q_J^* + t q_{J-1}^*)^2 s - (q_J + t q_{J-1})^2 s^*
&=\left((q_J^* + t q_{J-1}^*) v - (q_J + t q_{J-1}) v^*\right) \\ 
& \quad \cdot \left((q_J^* + t q_{J-1}^*) v + (q_J + t q_{J-1}) v^*\right).
\end{aligned}
\]
Now \eqref{system_one} and
\eqref{lb_v} give
\begin{equation} \label{system_two}
\begin{cases}
&|\, q_{J-1}^* v - q_{J-1} v^*|\ll Q^{2+8\phi-\delta/2+2\varepsilon},  \\
&\left|(q_J^* + t q_{J-1}^*) v - (q_J + t q_{J-1}) v^*\right| \ll Q^{2+8\phi-\delta/2+2\varepsilon}.
\end{cases}
\end{equation}
Multiplying the first inequality by $t$ and subtracting the second one yields
\[
\begin{cases}
&|q_{J-1}^* v - q_{J-1} v^*|\ll Q^{2+8\phi-\delta/2+2\varepsilon},\\
&\left|q_J^* v - q_J v^*\right| \ll Q^{2+8\phi-\delta/2+2\varepsilon}.
\end{cases}
\]

Consequently, there exist $X, Y \in \bR$ with 
\[
|X|,|Y| \ll Q^{2+8\phi-\delta/2+2\varepsilon}
\]
such that
\begin{equation} \label{s_equals_M}
\begin{pmatrix}
q_{J-1}^* & - q_{J-1}\\
q_J^* & - q_J
\end{pmatrix} \begin{pmatrix}
v\\v^*
\end{pmatrix}
= \begin{pmatrix}
X\\Y
\end{pmatrix}.
\end{equation}
In light of \eqref{lb_v} and~\eqref{s_equals_M}, Cramer's rule gives
\begin{align*} 
q_Jq_J^*\left|\frac{q_{J-1}}{q_J}-\frac{q_{J-1}^*}{q_J^*}\right| &= \left|\mathrm{det}\begin{pmatrix}
q_{J-1}^* & - q_{J-1} \\
q_J^* & - q_J\end{pmatrix}\right| \nonumber \\ &\ll Q^{1+2+8\phi-\delta/2+4\varepsilon-(1+\delta/2)}=Q^{2+8\phi-\delta+4\varepsilon}.
\end{align*}
By \eqref{denominator_bounds}, we thus have
\[ 
\left|\frac{q_{J-1}}{q_J}-\frac{q_{J-1}^*}{q_J^*}\right|\ll Q^{8\phi-\delta+8\varepsilon}. 
\]
It now follows from Proposition~\ref{prop_close_ratios}
that 
\begin{align*}
(\lambda^K\times \lambda^K)
(\cS_3) \ll Q^{(8\phi-\delta+8\varepsilon)(2\kappa-3-4(2-\kap)\varepsilon)+(1+2\varepsilon)c/(R+c)},  
\end{align*}
where the constant $c>0$ is as in Lemma~\ref{vq}. Coupling this with the trivial bound $|\fJ| \le 2,$ we find that 
\begin{equation}
\label{Contribution2223}
\int_{\cS_3} \fJ \dd \lambda^{2K} \ll Q^{(8\phi-\delta+8\varepsilon)(2\kappa-3-4(2-\kap)\varepsilon)+(1+2\varepsilon)c/(R+c)}. 
\end{equation}

\subsubsection{The size of $m_2$} 

Combining the various contributions \eqref{Contribution1}, \eqref{Contribution21}, \eqref{Contribution220}, \eqref{Contribution221},
\eqref{Contribution2221},
\eqref{Contribution2222},
\eqref{Contribution2223}
to $m_2$, we find that 
\[  
m_2 \ll  Q^{(2\eps - \phi)
(2 \kap - 3 - 4(2-\kap)\eps)
+ (1 + 2\eps) 2c / (R+c)}.  
\]

\subsection{Fourier decay} 

We are now in position to use the comparison Lemma \ref{ComparisonLemma} to deduce an upper bound for $\widehat{\nu}.$  Recall that for $k_1, k_2 \in \mathbb{Z}$ we have 
\[ \widehat{\nu}(k_1,k_2) = \int F(t,u)\d \nu(t,u)  
\]
where $F(t,u)$ is defined in \eqref{F_definition}. We showed in \S \ref{partial_d} that 
\[  
\max |\partial_t F| \ll \frac{|k_1|}{Q^{2-4\varepsilon}} + \frac{|k_2|}{Q^{1-2\varepsilon}} \ll \frac{|k_1|}{Q^{2-4\varepsilon}} + \frac{|k_2|Q}{Q^{2-2\varepsilon}} \ll Q^{\delta + 4\varepsilon} 
\] 
and  
\[    \max |\partial_u F| \ll \frac{|k_2|}{Q^{1-2\varepsilon}} =  \frac{|k_2|Q}{Q^{2-2\varepsilon}} \ll  Q^{\delta+2\varepsilon}  \ll Q^{\delta + 4 \varepsilon},
\] 
so we may apply Lemma \ref{ComparisonLemma} with 
\[ 
M_1 = M_2 = 
Q^{\delta + 4\varepsilon}.  
\]
With this choice of $M_1, M_2$, for the quantity $\Lambda(r)$ defined in Lemma~\ref{ComparisonLemma}, applying Lemma \ref{lem_FD} gives
\[ 
\Lambda(r) \ll \Big( \frac{r}{M_1}\Big)^{2\kappa -2 -4\varepsilon} = \frac{r^{2\kappa-2-4\varepsilon}}{Q^{(\delta+4\varepsilon)(2\kappa-2 -4\varepsilon)}}.
\]
Moreover, we demonstrated in \S \ref{m2calc} that the second moment 
$m_2$ of the function $F$ satisfies 
\[ 
m_2 \ll Q^{-\eta},
\]
where
\[
\eta = (\phi - 2\eps)
(2 \kap - 3 - 4\eps)
- (1 + 2\eps) 2c / (R+c).
\]

Combining these estimates, Lemma \ref{ComparisonLemma} gives
\begin{align*}  \widehat{\nu}(k_1, k_2) &\ll 
r + \frac{r^{2\kappa-2-4\varepsilon}}{Q^{(\delta+4\varepsilon)(2\kappa-2 -4\varepsilon)}}\Big(1 + \frac{Q^{\delta + 4\varepsilon}}{r} + \frac{Q^{2(\delta+4\eps) - \eta}}{r^4} \Big) \\ 
& = r + \frac{r^{2\kappa-2-4\varepsilon}}{Q^{(\delta+4\varepsilon)(2\kappa-2 -4\varepsilon)}} + \frac{r^{2\kappa - 3 -4\varepsilon}}{Q^{(\delta+4\varepsilon)(2\kappa -3-4\varepsilon) }} \\
& \qquad + 
\frac{r^{2\kappa - 6 -4 \eps}}
{Q^{(\delta + 4 \eps)(2\kappa -4 - 4\eps)
	+\eta}}. 
\end{align*}
We assume in the sequel that 
\[
\kap > 8/5,
\]
and we recall from \eqref{R_is_large} that 
$c < \eps R$.
Inserting the value of $\eta$ and simplifying gives
\begin{align}
\notag
\widehat{\nu}(k_1,k_2) &\ll  
r + \frac{r^{2\kappa-2-4\varepsilon}}{Q^{(\delta+4\varepsilon)(2\kappa-2 -4\varepsilon)}} + \frac{r^{2\kappa - 3 -4\varepsilon}}{Q^{(\delta+4\varepsilon)(2\kappa -3-4\varepsilon)}} \\ \notag 
& \qquad + 
\frac{r^{2\kappa - 6 -4 \eps}}
{Q^{(\delta + \phi + 2 \eps)(2\kappa - 3 - 4\eps) - (\del + 4 \eps) - (1+2\eps)2c/(R+c)}} 
\\ \notag
&\le
r + \frac{r^{2\kappa-2-4\varepsilon}}
{Q^{\del(2\kap-2 -4\eps)}} + \frac{r^{2\kappa - 3 -4\varepsilon}}{Q^{\del(2\kap -3-4\eps)}} + 
\frac{r^{2\kappa - 6 -4 \eps}}
{Q^{(\delta + \phi)
	(2\kappa - 3 - 4\eps) - 
	\del - 7 \eps}} \\
&\le
r + (rQ^{-\del})^{2\kap-3-4\eps} + 
\frac{r^{2\kappa - 6 -4 \eps}}
{Q^{(\delta + \phi)
	(2\kappa - 3 - 4\eps) - 
	\del - 7 \eps}}
.\label{lemma_gives} 
\end{align}
Setting
\[  
r = Q^{((\delta + \phi)
(2\kappa - 3 - 4\eps) - 
\del - 7 \eps)/
(2 \kap - 7 - 4 \eps)},
\] 
so that the first and last terms in \eqref{lemma_gives} are equal, we obtain 
\begin{align}
\label{obtain}
\widehat{\nu}(k_1, k_2) \ll  
Q^{- \: \frac{\del(22\kap-43-4\eps)-70\eps}
{10(7 - 2 \kap + 4 \eps)}} 
\le Q^{2 \eps - \: \frac{\del(22\kap-43-4\eps)}
{10(7 - 2 \kap + 4 \eps)}}.
\end{align}

\bigskip

We write 
\[ 
\|\bk\|_{\infty} = 
\max\{  |k_1|, |k_2| \}  
\]
for the supremum norm of the integer vector 
$\bk = (k_1, k_2).$ 
With \eqref{kk_size} in mind, we now distinguish two cases according to the relative sizes of 
$|k_1|$ and $|k_2|.$

\bigskip

\textbf{Case I.} The integers $k_1, k_2$ are such that $|k_1| \leq |k_2|.$ Then, by \eqref{kk_size}, we have
\begin{equation*} 
\max\{ |k_1|, |k_2|Q\} = |k_2| Q \asymp Q^{2+\delta} \end{equation*} 
and therefore $\|\bk\|_\infty = |k_2| \asymp Q^{1+\delta}.$  Now \eqref{obtain} gives 
\[ 
\widehat{\nu}(k_1, k_2) \ll \|\bk \|_\infty^{\frac{-\delta (22\kappa -43-4\eps)}{10(1+\delta)(7-2\kappa+4\eps)}+2\varepsilon}.  
\]

\bigskip

\textbf{Case II.} The integers $k_1, k_2$ satisfy $|k_1| > |k_2|.$ Then 
\[ \|\bk \|_\infty = |k_1| \le \max\{ |k_1|, |k_2|Q  \}  \asymp Q^{2+\delta},  \]
and combining with \eqref{obtain} we get
\begin{align*}
\widehat{\nu}(k_1, k_2) 
&\ll \| \bk \|_\infty^{\frac{-\delta (22\kappa -43)}
{10(2+\delta)(7-2\kappa)} + 4\varepsilon}. 
\end{align*}
Combining the two cases and using the fact that all norms on $\bR^2$ are equivalent, we deduce that 
\begin{equation} \label{nu_hat_ub}
\widehat{\nu}(k_1, k_2) \ll \|\bk \|_\infty^{\frac{-\delta (22\kappa -43-4\eps)}{10(2+\delta)(7-2\kappa+4\eps)}+2\varepsilon},  
\end{equation}
as required.

\begin{remark} \label{rem_justification_concrete_eta}
The upper bound~\eqref{nu_hat_ub} proves Theorem~\ref{MainThm}~\eqref{MainThm_point_decay}. One can see that if $\kappa$ is sufficiently close to $2$, and $\varepsilon>0$ is sufficiently small, and $\delta>0$ sufficiently large, then the exponent in the right-hand side of~\eqref{nu_hat_ub} can be made arbitrarily close to $1/30$. This justifies Remark~\ref{rem_concrete_eta}.
\end{remark}

\section{Lacunary approximation}
\label{LacSection}

In this section, we prove Theorem~\ref{GeneralCountingThm}. We closely follow~\cite{PVZZ2022}, adapting the calculations to our two-dimensional setup. 

We remark that, while in this article we only consider lacunary sequences, the same analysis can be done for smooth and, more generally, for $\alp$-separated sequences, cf. \cite{PVZZ2022}. 

Theorem~\ref{GeneralCountingThm} follows quite straightforwardly from Lemma~\ref{ebc}, Proposition~\ref{indieMainThm}, and Lemma \ref{lem_sum_mu_good} below. We refer the reader to~\cite[\S 3.3]{PVZZ2022} for further details.

Below, Lemma~\ref{ebc} recalls a classical tool for establishing asymptotic counting results. Then, we proceed to state and prove Proposition \ref{indieMainThm} and Lemma \ref{lem_sum_mu_good}, developing on the way the machinery needed.

Our proof is based on the following abstract counting result~\cite[Lemma~1.5]{harman}, which is a quantitative version of the divergence Borel--Cantelli lemma \cite{BV2023}.

\begin{lemma}\label{ebc}
Let $(X,\ca{B},\mu)$ be a probability space, let $(f_n(x))_{n \in \N}$ be a sequence of non-negative $\mu$-measurable functions defined on $X$, and let $(f_n)_{n \in \N },\ (\phi_n)_{n  \in \N}$ be sequences of real numbers  such that
\[
0\leq f_n \leq \phi_n \qquad (n \in \bN).
\]
Write
$$
\Phi(N)= \sum \limits_{n=1}^{N}\phi_n,
$$
and suppose that $\Phi(N) \to \infty$ as $N \to \infty$.
Suppose that, for arbitrary  $a,b \in \N$ with $a <  b$, we have
\begin{equation} \label{ebc_condition1}
\int_{X} \left(\sum_{n=a}^{b} \big( f_n(x) -  f_n \big) \right)^2\mathrm{d}\mu(x)\, \leq\,  C\!\sum_{n=a}^{b}\phi_n
\end{equation}

\noindent for an absolute constant $C>0$. Then, for any $\varepsilon>0$,  we have
\begin{equation} \label{ebc_conclusion}
\sum_{n=1}^N f_n(x)\, =\, \sum_{n=1}^{N}f_n\, +\, O\left(\Phi(N)^{\frac12}(\log(2+\Phi(N)))^{\frac{3}{2}+\varepsilon}+\max_{1\leq k\leq N}f_k\right)
\end{equation}
for $\mu$-almost all $x\in X$. %where $\Phi(N) = \sum\limits_{n=1}^{N}\phi_n$.
\end{lemma}

Lemma~\ref{ebc} allows to establish counting results in many situations. The key point in order to apply it is to verify hypothesis~\eqref{ebc_condition1}. This hypothesis is verified, in our situation of interest, by Proposition~\ref{indieMainThm} below. In order
to state it, we set
\begin{equation} \label{Et_def}
\cE_t = \{ (\alp, \gam) \in [0,1)^2: \| n_t \alp - \gam \| \le \psi(n_t) \}
\qquad (t \in \bN).
\end{equation}

%and~\ref{indiethm_general_alpha_separated} and Corollary~\ref{cor_indie_smooth}
%below allows to establish hypothesis~\eqref{ebc_condition1} for some more specific sequences of denominators $\cA$.

\begin{prop}  \label{indieMainThm} Let $\nu$ be a measure supported on $[0,1)^2$ satisfying~\eqref{LogDecay} 
%~\eqref{decaylac}
with $A>2$, and let $(n_t)_{t\in \bN} $ be a lacunary sequence of positive integers with $n_1 > 4$.
Let $\psi:\mathbb{N}\to[0,1/2)$. Then, for arbitrary  $a,b \in \mathbb{N}$ with $a < b$, we have 
\begin{equation}
\label{delta_set_frac23}
\begin{aligned}
	& 2 
	\mathop{\sum\sum}_{a\leq s<t\leq b} \nu(\cE_{s}\cap \cE_{t}) \le \left( \sum_{t=a}^{b} \nu(\cE_{t}) \right)^2 \\
	&+ O\left(\left(
	\sum_{t=a}^{b}\psi(n_t)\right)^{4/3}
	\log^+\left(\sum_{n=a}^{b}
	\psi(n_t)\right)
	+\sum_{t=a}^{b}\psi(n_t)\right).
\end{aligned} \end{equation}
\end{prop}

In order to prove %such results, 
Proposition~\ref{indieMainThm}, we approximate the left-hand side of \eqref{delta_set_frac23} by integrals of suitable functions; this is done in~\eqref{ub_mu_first}. Then, by employing Parseval's equality, we reduce the problem to bounding certain sums of Fourier coefficients, which is carried out in Proposition ~\ref{prop_sum_S_m_n}.

Given $\eps \in (0,1)$ and 
$\del \in (0,1/4)$,  we define  
$\chi_{\delta}: [0,1) \to \bR$ 
by   
\[
\chi_{\delta}(x):= \begin{cases}1 
\ \text{ if }  \ \|x\|\leq \delta  \\[1ex]
0 \ \text{ if } \  \|x\|>\delta,  \end{cases}   
\]
and its continuous upper and lower approximations 
$\chi_{\delta, \varepsilon}^\pm 
: [0,1)\to  \bR$  
by
\[
\chi_{\delta, \varepsilon}^+(x):
= \begin{cases} 1 & \text{\  if   }  \  \|x\|\leq \delta \\[1ex]
1+ \dfrac{1}{\delta\varepsilon}(\delta-\|x\|) &\text{\   if }  \  \delta < \|x\| \leq (1+\varepsilon)\delta \\[1.3ex]
0 & \text{\ if } \  \|x\|>(1+\varepsilon)\delta \end{cases}
\]
and
$$   
\chi_{\delta, \varepsilon}^-(x)
= \begin{cases} 1 & \text{\  if } \  \|x\|\leq (1 -\varepsilon)\delta  \\[1ex]
\dfrac{1}{\delta\varepsilon}(\delta-\|x\|) &\text{\   if } \  (1-\varepsilon)\delta < \|x\| \leq \delta \\[1.3ex]
0  & \text{\ if } \ \|x\|>\delta. \end{cases} $$
We now introduce the main difference compared to \cite{PVZZ2022}. Since therein $\gamma$ is fixed, whereas in the present article it is not, we require two-dimensional analogues of the functions employed in~\cite{PVZZ2022}. Given a real positive function $\psi:\mathbb{N}\rightarrow [0,1)$ and an integer $q\geq 4$, we define
$W_{q,\varepsilon}^{+}, W_{q,\varepsilon}^{-}: [0,1)^2\to\bR$ by \vspace{2mm}
\begin{equation*} %\label{Wdef}
W_{q,\varepsilon}^{+}(\alp,\gamma) :=\ %\Big( \sum_{p=0}^{q-1}\delta_{\frac{p+\gamma}{q}} (x)  \Big) *\chi_{\frac{\psi(q)}{q},\varepsilon}^{+}(x)=
\sum_{p=0}^{q-1}\chi_{\frac{\psi(q)}{q}, \varepsilon}^{+}\left(\textstyle{\alp- \frac{p+\gamma}{q}}\right)
\end{equation*}
and
\begin{equation*} 
W_{q,\varepsilon}^{-}(\alp,\gamma)  :=\sum_{p=0}^{q-1}\chi_{\frac{\psi(q)}{q}, \varepsilon}^{-}\left(\textstyle{\alp- \frac{p+\gamma}{q}}\right)   \, . \vspace{2mm}
\end{equation*}
In the sequel, we write 
$ W_{q,\varepsilon}^{\pm}$ to mean both the `upper' and `lower' functions, when this creates no  risk of confusion. Similarly we  write $\chi_{\delta, \varepsilon}^{\pm}$ to refer to both $\chi_{\delta, \varepsilon}^{+}$ and $\chi_{\delta, \varepsilon}^{-}$. As $\psi$ is considered fixed, we write $W_{q,\varepsilon}^{\pm}$ instead of the more descriptive $W_{q,\varepsilon, \psi}^{\pm}$.

It follows that for any $0<\varepsilon \le 1$ and any integer $t\in\bN$, \vspace{2mm}
\begin{equation} \label{muineq}
\int_{[0,1)^2}W_{n_t,\varepsilon}^-(\alp,\gamma)\mathrm{d}\nu(\alp,\gamma)
\; \leq  \;
\nu(\cE_t) \;  \leq \; \int_{[0,1)^2}W_{n_t,\varepsilon}^+(\alp,\gamma)\mathrm{d}\nu(\alp,\gamma). \vspace{2mm}
\end{equation}
As in~\cite{PVZZ2022}, we evaluate these integrals using the Fourier series expansions of $W_{q,\varepsilon}^{\pm}$. For 
$\mathbf{k} = (k_1,k_2) \in \bZ$ let $ \widehat{\chi}_{\delta, \varepsilon}^{\pm}(\mathbf{k})$ and $\widehat{W}_{q,\varepsilon}^{\pm}(\mathbf{k})$ denote the $\mathbf{k}^{\mathrm{th}}$ Fourier coefficient of $ {\chi}_{\delta, \varepsilon}^{\pm} $ and $ {W}_{q,\varepsilon}^{\pm},$
respectively.  

We denote
\[
A(q,k_1):=\dfrac{q\left(\cos(2\pi k_1\psi(q)q^{-1})-\cos(2\pi k_1\psi(q)q^{-1}(1+\varepsilon))\right)}{2\pi^2 k_1^2\psi(q)q^{-1}\varepsilon}.
\]
Using \cite[Equation (46)]{PVZZ2022} we compute, for $k_1\neq 0$, that
\begin{equation} 
\label{fcoef}
\begin{aligned}
&\widehat{W}_{q,\varepsilon}^+(k_1,k_2)\, =
\int_0^1\int_0^1 W_{q,\varepsilon}^{+}(\alp,\gamma) e\left(-k_1 \alp
-k_2\gamma)\right) \dd \alp \dd\gamma
\\&=
\,\begin{cases}
	\int_0^1 e\left( \dfrac{-k_1 \gamma}{ q}- k_2\gamma\right) A(q,k_1)\dd\gamma, &\text{if } q \mid k_1 \\[3ex]
	0, &\text{if } q\nmid k_1
\end{cases} \\&=
\,\begin{cases}
	%\dfrac{q\left(\cos(2\pi k_1\psi(q)q^{-1})-\cos(2\pi k_1\psi(q)q^{-1}(1+\varepsilon))\right)}{2\pi^2 k_1^2\psi(q)q^{-1}\varepsilon}
	A(q,k_1), &\text{if }
	k_1=-q k_2 \\[3ex]
	0, &\text{if } k_1\neq - q k_2.
\end{cases}
\end{aligned}
\end{equation}
We also find that
\begin{equation} \label{fcoef_zero}
\widehat{W}_{q,\varepsilon}^{+}(0,0)=(2+\varepsilon) \psi(q)
\end{equation}
and
\begin{equation} \label{fcoef_zero_nonzero}
\widehat{W}_{q,\varepsilon}^{+}(0,k_2)=0
\qquad (k_2 \ne 0).
\end{equation}

A similar analysis applies to $\widehat{W}_{q,\eps}^-$. We denote
\[
B(q,k_1):=\dfrac{q\left(\cos(2\pi k_1\psi(q)q^{-1}(1-\varepsilon))-\cos(2\pi k_1\psi(q)q^{-1}) \right)}{2\pi^2 k_1^2\psi(q)q^{-1}\varepsilon}.
\]
For $k_1\neq 0$, we have
\begin{equation} \label{fcoef2}
\widehat{W}_{q,\varepsilon}^-(k_1,k_2)\, =\,
\begin{cases}
B(q,k_1)  &\text{\  if   }  k_1= - q k_2 \\[3ex]
0 & \text{\  if   }  k_1\neq - q k_2.
\end{cases}
\end{equation}
We also find that
\begin{equation} \label{fcoef_zero2}
\widehat{W}_{q,\varepsilon}^-(0,0)  = (2 -  \varepsilon)\, \psi(q)  
\end{equation}
and
\begin{equation} \label{fcoef_zero_nonzero2}
\widehat{W}_{q,\varepsilon}^-(0,k_2)  = 0  
\qquad (k_2 \ne 0).
\end{equation}

Now $\sum\limits_{\mathbf{k} \in \bZ^2}| \widehat{W}_{q,\varepsilon}^{\pm}(\mathbf{k} )  |  < \infty $, so the Fourier series
\[ \sum_{\mathbf{k} \in \bZ^2}\widehat{W}_{q,\varepsilon}^{\pm}(\mathbf{k}) e( k_1 \alp + k_2 \gamma ) \]
converges uniformly to $W_{q,\varepsilon}^{\pm}(\alp,\gamma)$ for all $(\alp,\gamma) \in [0,1)^2$. It  follows that
\begin{eqnarray*}
\int_0^1 W_{q,\varepsilon}^{\pm}(\alp,\gamma)  \; \mathrm{d}\nu(\alp,\gamma) \; =  \;  \sum_{\mathbf{k}\in \bZ^2} \, \widehat{W}_{q,\varepsilon}^{\pm}(\mathbf{k})  \;  \widehat{\nu}(-\mathbf{k})   \, .
\end{eqnarray*}
Together with \eqref{muineq}, \eqref{fcoef_zero},  \eqref{fcoef_zero2} and the fact that $\widehat{\nu}(0,0)=1$, this implies that if $t \in \bN$ then
\begin{equation} \label{mu_ie}
\begin{array}{ll}
\nu( \cE_t)  \  \leq \ (2+\varepsilon) \, \psi(n_t)  \  + \displaystyle{\sum_{\mathbf{k} \in \bZ^2 \setminus \{\mathbf{0} \} }}\widehat{W}_{n_t,\varepsilon}^{+}(\mathbf{k})  \; \widehat{\nu}(-\mathbf{k}),
\\[5ex]
\nu( \cE_t)  \ \geq  \  (2-\varepsilon) \,  \psi(n_t)   \ + \displaystyle{\sum_{\mathbf{k} \in \bZ^2 \setminus \{\mathbf{0} \} }}  \widehat{W}_{n_t,\varepsilon}^{-}(\mathbf{k}) \; \widehat{\nu}(-\mathbf{k}).
\end{array}
\end{equation}

The arguments leading to \cite[Equations (51) and (52)]{PVZZ2022} can be used to deduce from~\eqref{fcoef} and~\eqref{fcoef2} the inequalities
\begin{eqnarray} \label{ub_Wpm_zero}
|\widehat{W}_{q,\varepsilon}^{\pm}(qk_2,-k_2)| &\leq& (2+\varepsilon) \, \psi(q),\label{W_ub_psi}\\
\label{ub_Wpm_nonzero}
|\widehat{W}_{q,\varepsilon}^{\pm}(q k_2,-k_2)| &\leq& \frac{1}{\pi^2k_2^2\psi(q)\varepsilon} \label{W_ub_1s2},
\end{eqnarray}
valid for any integer $k_2\neq 0$.

\begin{lemma} \label{main_part_bounded}
Let $0< \varepsilon,\tilde{\varepsilon} \le 1$.  Then, for any integers $q, r \geq 4$, we have
\begin{equation} \label{ub_single_sum}
\sum_{s \in \bZ} \big|\widehat{W}_{q,\varepsilon}^{\pm}(sq,-s)\big| \  <  \ \frac{12}{\varepsilon^{1/2}}
\end{equation}
and
\begin{equation} \label{ub_double_sum}
\sum_{s\in \bZ}   \sum_{t\in \bZ}  \big|\widehat{W}_{q,\varepsilon}^{\pm}(sq,-s)
\widehat{W}_{r,\tilde{\varepsilon}}^{\pm}(tr,-r)\big|  \ \leq \  \frac{144}{\varepsilon^{1/2}\cdot\tilde{\varepsilon}^{1/2}}  \, .
\end{equation}
\end{lemma}
\begin{proof}
We have
\begin{align*}
&\sum_{s \in \bZ}
\big|\widehat{W}_{q,\varepsilon}^{\pm}(sq,-s)\big| \\
&=
\sum_{|s| \le \frac{1}{\sqrt{2}\pi\psi(q)\varepsilon^{1/2}}}  \big|\widehat{W}_{q,\varepsilon}^{\pm}(sq,-s)\big| +
\sum_{|s| > \frac{1}{\sqrt{2}\pi\psi(q)
		\varepsilon^{1/2}}}     
\big|\widehat{W}_{q,\varepsilon}^{\pm}(sq,-s)\big|.
\end{align*}
Using \eqref{fcoef_zero} and \eqref{W_ub_psi} in the first sum and \eqref{W_ub_1s2}
in the second, we obtain
\begin{equation*} 
\begin{aligned}
	\sum_{s \in \bZ}
	\big|\widehat{W}_{q,\varepsilon}^{\pm}(sq,-s)\big| &\leq    
	\: 3 
	\sum_{|s| \le \frac{1}{\sqrt{2}\pi\psi(q)
			\varepsilon^{1/2}}}   \hspace{-5mm} \psi(q)  \ \ + 
	\hspace*{3ex}  
	\sum_{|s| >
		\frac{1}{\sqrt{2}\pi\psi(q)
			\varepsilon^{1/2}}} 
	\frac{1}{\pi^2 s^2\psi(q)\varepsilon} \\
	&\leq \:
	\frac{3 \sqrt{2}}{\pi\varepsilon^{1/2}}  \ + 3\psi(q) + \ \frac{ 2 \sqrt{2}}{\pi\varepsilon^{1/2}} + 6\psi(q) \ <  \  \frac{12}{\varepsilon^{1/2}}.
\end{aligned}
\end{equation*}
The inequality \eqref{ub_double_sum} follows from \eqref{ub_single_sum} and the fact that 
\begin{align*}
&\sum_{s\in \bZ} \sum_{t\in \bZ}  \big|\widehat{W}_{q,\varepsilon}^{\pm}(sq,-s)\widehat{W}_{r,\tilde{\varepsilon}}^{\pm}(tr,-r)\big|
\\  &=
\Big(\sum_{s \in \bZ} \big|\widehat{W}_{q,\varepsilon}^{\pm}(sq,-s)\big|\Big)
\Big(\sum_{t  \in \bZ} \big|\widehat{W}_{r,\tilde{\varepsilon}}^{\pm}(tr,-r)\big|\Big).
\end{align*}
\end{proof}

\begin{lemma} \label{lem2}
Let $\nu$ be a probability measure supported on a subset $F$ of  $\I^2$. Let  $\cE_{t}$ be given by~\eqref{Et_def}. Then, assuming that $n_t\geq 4$, we have
\begin{eqnarray}
\nu(\cE_{t}) &\leq& 3\psi(n_t) + 12 \max_{k\in\bZ\setminus\{0\}}|\widehat{\nu}(-k n_t,k)| ,
\label{lem2_result1}
\end{eqnarray}
and
\begin{eqnarray}
\nu(\cE_{t}) &\leq& 3\psi(n_t) + 2\sum_{k=1}^{\infty}\frac{|\widehat{\nu}(-k n_t,k)|}{k}.  \label{lem2_result2}
\end{eqnarray}
\end{lemma}

\begin{proof}
It follows from \eqref{mu_ie} with $\varepsilon= 1$, and \eqref{fcoef}, that
\begin{equation} \label{lem2_ie2}
\nu(\cE_{t})   \ \le  \  3\psi(n_t)     + \sum_{k \in \bZ \setminus \{0 \}  }\widehat{W}_{n_t,1}^{+}(k n_t, -k)   \  \widehat{\nu}(-k n_t, k).
\end{equation}
The desired estimate~\eqref{lem2_result1} now immediately follows from this, together with the fact that by  Lemma~\ref{main_part_bounded},
\begin{equation*} \label{lem2_ie3}
\Big|\sum_{k \in \bZ \setminus \{0 \}  }\widehat{W}_{n_t,1}^{+}(k n_t, -k)   \  \widehat{\nu}(-k n_t,k)\Big|   \  \stackrel{\eqref{ub_single_sum}}{\leq}  \ 12 \, \max_{k\in\bZ\setminus \{0\} }|\widehat{\nu}(-k n_t,k)|.
\end{equation*}
To prove \eqref{lem2_result2}, note that the geometric mean of the inequalities~\eqref{W_ub_psi} and~\eqref{W_ub_1s2} implies that, for every $k\in\bZ\setminus\{0\}$,
$$
\big|\widehat{W}_{n_t,1}^{\pm}(k n_t,-k)\big|  \ \leq \  \frac{\sqrt{3}}{\pi |k|}    \ \leq \   \frac{1}{|k|}.
$$
This and \eqref{lem2_ie2} imply \eqref{lem2_result2}.
\end{proof}

The following convergence result readily follows from Lemma~\ref{lem2}.

\begin{thm} \label{mainCONV}  Let $\nu$ be a probability  measure supported on a subset $F$ of
$ \ \I^2 \, $.  Let $ (n_t)_{t\in \N} $ be a sequence of natural numbers, and let $\psi:\mathbb{N}\rightarrow \I$. Suppose that at least one of the following two conditions is satisfied:
\begin{equation} \label{cond2_2}
\sum_{t=1}^{\infty} \; \max_{k\in\bZ \setminus \{ 0\} }|\widehat{\nu}(-k n_t,k)|< \infty,
\end{equation}
\begin{equation} \label{cond2}
\sum_{t=1}^{\infty}   \, \sum_{k\in\bZ \setminus
	\{ 0\} }  \frac{|\widehat{\nu}
	(-kn_t,k)|}{|k|}
< \infty .
\end{equation}
Then
\begin{equation*} %\label{convcondition}
\nu(\cW(\psi))=0  
\end{equation*}
if
\begin{equation} \label{convergence}
\sum\limits_{n=1}^{\infty}\psi(n_t)< \infty.
\end{equation}
\end{thm}
\begin{proof}
It follows from Lemma~\ref{lem2} and hypotheses~\eqref{cond2_2} and~\eqref{cond2} of the theorem that
\[
\sum\limits_{t=1}^{\infty}\nu(\cE_t)<\infty.
\]
The result now follows from the first Borel--Cantelli lemma.
\end{proof}

\begin{remark}
Note that the condition \eqref{LogDecay} implies \eqref{cond2_2}, as well as \eqref{cond2}. Thus, if we have~\eqref{convergence}, then it follows from Theorem~\ref{mainCONV} that
\[
\cN(\alp,\gam,\psi,N) \ll 1,
\]
for $\nu$-almost all  $(\alpha,\gamma)$. At the same time, under 
\eqref{convergence},
the right-hand side of~\eqref{countlacresult} is also bounded. We thereby obtain the conclusion of Theorem
\ref{GeneralCountingThm} under the assumption \eqref{convergence}.
We may therefore assume for the remainder of this section that
\begin{equation*} \label{divergence}
\sum_{t=1}^\infty \psi(n_t) = \infty.
\end{equation*}
\end{remark}

\begin{cor} \label{cor_min_psi}
We work within the framework of Theorem~\ref{thm_lacunary}. Let $\tau>1$, let $\psi_\tau: \bN \to [0,1)$ be any function such that
\[
\psi_\tau(n_t) = t^{-\tau}
\qquad (t \in \bN),
\]
and put
$\widetilde{\psi} = 
\max \{ \psi, \psi_\tau \}$.
Then, for $\nu$-almost all pairs $(\alp, \gam)$, the formula \eqref{countlacresult} for $\psi$ is equivalent to the same formula with $\psi$ replaced by $\widetilde{\psi}$. In particular, when proving Theorem \ref{thm_lacunary}, we may assume that
\begin{equation*}
\label{psi_is_not_small}
\psi(n_t) \geq t^{-\tau}
\qquad (t \in \bN).
\end{equation*}
\end{cor}

\begin{proof}
It follows from the definition~\eqref{def_cN} of the function $\cN(\alp,\gam,\psi,N)$ that
\[
\cN(\alp,\gam,\psi,N)\leq \cN(\alp,\gam,\widetilde{\psi},N)\leq \cN(\alp,\gam,\psi,N)+\cN(\alp,\gam,\psi_{\tau},N)
\]
At the same time, it follows from~Theorem~\ref{mainCONV} that, for $\nu$-almost all pairs $(\alp,\gam)$, we have $(\alp, \gam) \notin \cW(\psi_{\tau})$, which is equivalent to
\[
\cN(\alp,\gam,\psi_{\tau},N)
\ll 1,
\]
Thus, $\nu$-almost surely,
\[
\cN(\alp,\gam,\widetilde{\psi},N)=\cN(\alp,\gam,\psi_{\tau},N) + O(1).
\]
For these pairs $(\alpha, \gam)$, when we change $\psi$ to $\widetilde{\psi}$, we modify the left-hand side of \eqref{countlacresult} by $O(1)$.
Similarly, using the notation $\widetilde{\Psi}:=\sum_{t \le N} \widetilde{\psi}(n_t)$, we have
\[
\Psi(N)\leq\widetilde{\Psi}(N)\leq\Psi(N)+\sum_{t=1}^N\psi_{\tau}(n_t)=\Psi(N)+O(1),
\]
so by changing $\psi$ to $\widetilde{\psi}$ we modify the right-hand side of \eqref{countlacresult} by $O(1)$.
The term $O(1)$ is controlled by the error term in the right-hand side of \eqref{countlacresult}, so the claim is proved.
\end{proof}

\begin{lemma} \label{lem_sum_mu_good}
Let $\nu$ be a probability  measure supported on a subset $F$ of
$ [0,1)^2 $.  Let $  (n_t)_{t\in \bN} $ be a lacunary sequence of natural numbers such that $n_1 > 4$, and let $\psi:\mathbb{N}\rightarrow [0,1/2)$.
Suppose there exists a constant $A$ so that the measure $\nu$ satisfies \eqref{LogDecay}.
Then, for arbitrary $a,b \in \bN$ with $a < b$, we have
\begin{equation} \label{lem_sum_mu_good_result}
\sum_{t=a}^b\nu(\cE_t) = 2\sum_{t=a}^b\psi(n_t) \, + \, O\left( \min\Big(1,\sum_{t=a}^b\psi(n_t)\Big) \right).
\end{equation}
\end{lemma}

\begin{proof}
Let $\eps_a, \ldots, \eps_b \in (0,1].$
Combining \eqref{mu_ie} with \eqref{fcoef} and \eqref{fcoef2} yields
\begin{multline*}
\left| \sum_{t=a}^{b}\nu(\cE_t)  - 2\sum_{t=a}^b \psi(n_t) \right| \,\, \ls \,\, \sum_{t=a}^b \psi(n_t)\varepsilon_t \\
+  \max_{\circ \in \{+,-\}}
\sum_{t=a}^b\Big| 
\sum_{k\neq 0}
\widehat{W}_
{n_t,\varepsilon_t}^{\circ}
(kn_t,-k)\widehat{\nu}
(-kn_t, k) \Big|.
\end{multline*}
By \eqref{LogDecay}, we have 
$| \widehat{\nu}(kn_t, -k)|\ll t^{-A},$ and combining this with \eqref{ub_single_sum} gives
\begin{align*} 
\left| \sum_{t=a}^{b}\nu(\cE_t)  - 2\sum_{t=a}^b \psi(n_t) \right|
\ls \sum_{t=a}^b \psi(n_t)
\varepsilon_t + \sum_{t=a}^b\frac{O(1)}{t^{A}\sqrt{\varepsilon_t}}.
\end{align*}
This inequality is essentially the same as \cite[Equation (69)]{PVZZ2022}, and the arguments given there establish the claim. 
\end{proof}

\bigskip

Now we have all of the tools that we need to prove Proposition \ref{indieMainThm}. Observe from the definition that, for $s,t\in\mathbb{N}$, we have
\begin{equation} \label{ub_mu_first}
\nu(\cE_{s}\cap \cE_{t}) \leq \int_{0}^{1} W_{n_t,\varepsilon_t}^+(\alp,\gamma) W_{n_s,\varepsilon_s}^+(\alp,\gamma) \, \mathrm{d}\nu(\alp,\gamma).
\end{equation}
We introduce the notation
\[
W_{s,t}^+(\alp,\gamma):=W_{n_t,\varepsilon_t}^+(\alp,\gamma) W_{n_s,\varepsilon_s}^+(\alp,\gamma),
\]
and then we can compute the right-hand side of~\eqref{ub_mu_first} as follows:
\begin{align} 
\notag
&\int_{0}^{1} W_{n_t,\varepsilon_t}^+(\alp,\gamma) W_{n_s,\varepsilon_s}^+(\alp,\gamma) \, \mathrm{d}\nu(\alp,\gamma)     
\\ \notag & \qquad =
\int_{0}^{1} W_{s,t}^+(\alp,\gamma) \, \mathrm{d}\nu(\alp,\gamma)
\\ 
\label{Wst_Fourier} & \qquad =
\widehat{W}_{s,t}^+(0,0) + \sum_{\mathbf{k}\in\bZ^2\setminus\{\mathbf{0}\}}\widehat{W}^+_{s,t}(\mathbf{k})\widehat{\nu}(-\mathbf{k}).
\end{align}
Writing
\[
\lam \cE_t = \{
(\alp, \gam) \in [0,1)^2:
\| n_t \alp - \gam \| \le \lam \psi(n_t) \}
\]
for $\lam > 1$ and $t \in \bN$,
\begin{equation} \label{main_q_Harman_ub2}
\begin{aligned}
\widehat{W}_{s,t}^+(0,0)
&=\int_0^1\int_0^1 W_{n_s,\varepsilon_s}^+(\alp,\gamma)W_{n_t,\varepsilon_t}^+(\alp,\gamma) \dd \alp \dd \gamma.
\\&\leq\left| (1+\varepsilon_s)\cE_s \cap (1+\varepsilon_t)\cE_t \right|
\\&=
4(1+\varepsilon_s)(1+\varepsilon_t)\psi(n_s)\psi(n_t) 
\\ &\qquad + O\left( (n_s,n_t)\min\Big(\frac{\psi(n_s)}{n_s},\frac{\psi(n_t)}{n_t}\Big)  \right),
\end{aligned}
\end{equation}
where in the last step we have used \cite[Equation~(3.2.5)]{harman}. 

Further, we use the fact that for $\mathbf{k}\in\mathbb{Z}^2$, we have $\widehat{W}_{n_s,\varepsilon_s}^+(\mathbf{k})\neq 0$ only if $\mathbf{k}=(n_s k,-k)$ for some $k\in\bZ$. We also use the convolution theorem. Hence, for the second sum in~\eqref{Wst_Fourier}, 
\begin{multline} \label{Fourier_sum_principal}
\sum_{\mathbf{k}\in\bZ^2\setminus\{\mathbf{0}\}}\widehat{W}^+_{s,t}(\mathbf{k})\widehat{\nu}(-\mathbf{k})
=
\sum_{\substack{\mathbf{k},
	\mathbf{l}\in\bZ^2  \\ \mathbf{k}+\mathbf{l}\ne \mathbf{0}}}\widehat{W}_{n_s,\varepsilon_s}^+(\mathbf{k}) \widehat{W}_{n_t,\varepsilon_t}^+(\mathbf{l})\widehat{\nu}
(-\mathbf{k} - \mathbf{l})
\\ =
\sum_{\substack{k,\ell\in\bZ \\ k+\ell\ne 0\text{ or } \\ n_s k+n_t \ell \ne 0}}\widehat{W}_{n_s,\varepsilon_s}^+(n_s k, -k) \widehat{W}_{n_t,\varepsilon_t}^+(n_t \ell, -\ell)\widehat{\nu}(-n_s k - n_t \ell, k+\ell).
\end{multline}

\begin{prop}\label{prop_sum_S_m_n} We denote by $S_{s,t}$ the sum in the right-hand side of~\eqref{Fourier_sum_principal}.
Then, for any  $s,t \in \N$ with $s < t$,
\begin{equation} \label{S_bound}
|S_{s,t}| \ll
\frac{\psi(n_s)}{t^{A}\varepsilon_t^{1/2}}\, +\, 
\frac{\psi(n_t)}{s^{A}\varepsilon_s^{1/2}}
\, + \, \frac{1}{t^{A}\varepsilon_s^{1/2}
	\eps_t^{1/2}}   \,  +  \ |T(s,t)|,
\end{equation}
where 
\begin{align} \label{def_S_prime}
\notag &T(s,t) =  \\
&\mathop{ \mathop{\sum\sum}_{ k,\ell\in\mathbb{Z}\setminus \{0\} } }_{%1\leq 
	|kn_s+\ell n_t|<n_s^{1/2}} \!\!\! \what{W}_{n_s,\varepsilon_s}^{+}(kn_s,-k)\what{W}_{n_t,\varepsilon_t}^{+}(\ell n_t,-\ell)  \what{\nu}\left( -kn_s-\ell n_t,k+\ell\right).
\end{align}
\end{prop}

\begin{proof}
%In order to deduce the claimed upper bound, 
We decompose the sum $S_{s,t}$ as
\[
S_{s,t} = %S_0(s,t) +
S_1(s,t) + S_2(s,t) + S_3(s,t) + T(s,t),
\] where $T(s,t)$ is as above and we define 
\begin{align*}
&S_1(s,t) \\ &= \mathop{\sum\sum}_{\substack{ k,\ell\in\mathbb{Z} \\ k \ell = 0 \\ k n_s + \ell n_t \ne 0}}
\what{W}_{n_s,\varepsilon_s}^{+}(kn_s,-k)\what{W}_{n_t,\varepsilon_t}^{+}(\ell n_t,-\ell)  \what{\nu}\left( -kn_s-\ell n_t,k+\ell\right),\\
&S_2(s,t) \\ &= \!\!\!\!\! \mathop{ \mathop{\sum\sum}_{ k,\ell\in\mathbb{Z}\setminus \{0\} } }_{ |kn_s+\ell n_t|\gs n_t/2} \!\!\!\!\what{W}_{n_s,\varepsilon_s}^{+}(kn_s,-k)\what{W}_{n_t,\varepsilon_t}^{+}(\ell n_t,-\ell)  \what{\nu}\left( -kn_s-\ell n_t,k+\ell\right), \\
&S_3(s,t) \\ &= \mathop{ \mathop{\sum\sum}_{ k,\ell\in\mathbb{Z}\setminus \{0\} } }_{n_s^{1/2} \ls|kn_s+\ell n_t| < n_t/2} \hspace{-6mm}\what{W}_{n_s,\varepsilon_s}^{+}(kn_s,-k)\what{W}_{n_t,\varepsilon_t}^{+}(\ell n_t,-\ell)  \what{\nu}\left(-kn_s-\ell n_t,k+\ell\right).
\end{align*}

To bound $S_1(s,t)$ from above, we first consider the terms with $\ell=0$. By the
lacunarity of the sequence, and~\eqref{LogDecay},
\[
|\widehat{\nu}(-k n_s,k)|\ll \left(\log|kn_s|\right)^{-A}\ll |s|^{-A} \qquad (s \in \bN).
\]
Together with~\eqref{fcoef_zero} and \eqref{ub_single_sum}, this
implies that the sum $S_1(s,t)$ restricted to $\ell = 0$ is 
\[
\ll \  \frac{(2+\varepsilon_t)\psi(n_t)}{|s|^{A}}\sum_{k\in \bZ }  \widehat{W}_{n_s,\varepsilon_s}^{+}(k n_s,-k)  \ \ll \
\frac{\psi(n_t)}{|s|^{A}\varepsilon_s^{1/2}}\, .
\]
A similar bound holds for the terms with $k=0$, and adding up we deduce that
\begin{equation} \label{ub_S1}
S_1(s,t) \ll
\frac{\psi(n_t)}{|s|^{A}\varepsilon_s^{1/2}} + \frac{\psi(n_s)}{|t|^{A}\varepsilon_t^{1/2}}.
\end{equation}

For $S_2(s,t),$ the condition $|k n_s+\ell n_t|\geq n_t/2$, together with the lacunarity of the sequence,
and \eqref{LogDecay}, gives
\[
|\widehat{\nu}\left( -k n_s - \ell n_t,k+\ell\right)| \ll |t|^{-A}.
\]
This and \eqref{ub_double_sum} imply that
\begin{equation} \label{ub_S2}
S_2(s,t) \ll \frac{1}{|t|^{A}\varepsilon_s^{1/2}\varepsilon_t^{1/2}}.
\end{equation}

For $S_3(s,t),$ the hypothesis  $n_s^{1/2}\leq |k n_s+\ell n_t|< n_t/2$ together with the
lacunarity of $(n_t)_{t\in \bN}$
and~\eqref{LogDecay} gives
\[
|\widehat{\nu}\left( -k n_s- \ell n_t, k+\ell \right)| \ll |s|^{-A}.
\]
Observe that under the condition of summation we have
\[
|k n_s/n_t+\ell|< 1/2;
\]
this implies that  for each  $k$ there exists at most one index $\ell=\ell_k$ such that the pair $(k,\ell_k)$ satisfies the condition of summation. Moreover, the signs of $k$ and $\ell=\ell_k$ are necessarily opposite. These observations allow us to deduce the upper bound
\begin{equation*}
S_3(s,t) \ll \frac{1}{|s|^{A}} \mathop{\mathop{\sum}_{ k\in\bN \, : } }_{\ell_k {\rm \, exists}} |\widehat{W}_{n_s,\varepsilon_s}^{+}(k n_s,-k)| 
\cdot
|\widehat{W}_{n_t,\varepsilon_t}^{+}(\ell_k n_t,-\ell_k) |.
\end{equation*}
We proceed by applying~\eqref{ub_Wpm_zero} to $\widehat{W}_{n_t,\varepsilon_t}^{+}(\ell_k n_t,-\ell_k)$
and then the  bound \eqref{ub_single_sum} to the sum  $\sum\limits_{k\in\bN}|\widehat{W}_{n_s,\varepsilon_s}^{+}(k n_s,-k)|$ to eventually get \vspace{-2mm}
\begin{equation}  \label{ub_S3}
S_3(s,t) \ll|s|^{-A}\psi(n_t)\varepsilon_s^{-1/2}.
\end{equation}
Combining \eqref{ub_S1}, \eqref{ub_S2} and~\eqref{ub_S3} yields the result.
\end{proof}

\begin{prop}  \label{indiethm_general} 
Under the assumptions of Proposition~\ref{indieMainThm}, let $\tet \in (0,1]$. Then, for any  $a,b \in \mathbb{N}$ with $a <  b$, we have 
\begin{multline}
\label{ub_sum_eps_psipsi_combined}
2 \mathop{\sum\sum}_{a\leq s<t\leq b} \nu(\cE_{s}\cap \cE_{t})
\ls \left( \sum_{t=a}^{b}\psi(n_t)\right)^2 + O\left( 
\sum_{t=a}^b \psi(n_t)
\right)
\\ +  O\Bigg(\left(\sum_{t=a}^b\psi(n_t) \right)^{2-\theta}\!\!\!\log^+ \left(\sum_{t=a}^b\psi(n_t) \right)\,\, +\,\,  \mathop{\sum\sum}_{a\leq s<t\leq b}T(s,t) 
\\   +  \,\, \mathop{\sum\sum}_{a\leq s<t\leq b}(n_s,n_t)\min\Big(\frac{\psi(n_s)}{n_s},\frac{\psi(n_t)}{n_t}\Big) \Bigg).
\end{multline}
\end{prop}

\begin{proof}
By using the upper bound \eqref{ub_mu_first}, and then \eqref{main_q_Harman_ub2}, \eqref{S_bound} and Proposition~\ref{prop_sum_S_m_n}  with
\begin{equation} \label{def_varepsilon_t}
\varepsilon_t =
\min\Bigg(2^{-\theta},\, \Big( \sum_{s=a}^{t}\psi(n_s) \Big)^{-\theta}\Bigg)
\qquad (a \le t \le b),
\end{equation}
we find
\begin{multline} \label{ub_sum_mu}
2 \mathop{\sum\sum}_{a\leq s<t\leq b} \nu(\cE_{s}\cap \cE_{t}) \leq
\mathop{\sum\sum}_{a\leq s<t\leq b}
\Bigg(4(1+\varepsilon_s)(1+\varepsilon_t)\psi(n_s)\psi(n_t)+\frac{\psi(n_s)}{t^{A}\varepsilon_t^{1/2}}\,\\
+\, 
\frac{\psi(n_t)}{s^{A}\varepsilon_s^{1/2}}
\, + \, \frac{1}{t^{A}\varepsilon_s^{1/2}\varepsilon_t^{1/2}}+T(s,t)\Bigg)\\
+ O\left(\mathop{\sum\sum}_{a\leq s<t\leq b} (n_s,n_t)\min\Big(\frac{\psi(n_s)}{n_s},\frac{\psi(n_t)}{n_t}\Big)\right).
\end{multline}
It follows from \eqref{def_varepsilon_t} that
\begin{equation} \label{epsilon_inverse_bound}
\varepsilon_t^{-1} \leq \max\left(2^{\theta},t^{\theta}\right) < 2t.
\end{equation}
Consequently $\sum\limits_{t=1}^{\infty}1/(t^{A}\varepsilon_t^{1/2}) <\infty,$ so we have the  upper bound
\begin{equation} \label{ub_first_two_terms}
\sum_{a\leq s<t\leq b}\left(\frac{\psi(n_t)}{s^{A}\varepsilon_s^{1/2}} + \frac{\psi(n_s)}{t^{A}\varepsilon_t^{1/2}}\right)
\ll
\sum_{t=a}^b \psi(n_t)
\end{equation}
on part of the right-hand side of \eqref{ub_sum_mu}.

To bound the double sum $\sum\limits_{a\leq s<t\leq b} t^{-A}\varepsilon_s^{-1/2}\varepsilon_t^{-1/2}$ we consider two cases, according to whether $\sum\limits_{t=a}^{b}\psi(n_t)>2$ or $\sum\limits_{t=a}^{b}\psi(n_t)\leq 2$.

\noindent Case A: $\sum\limits_{t=a}^{b}\psi(n_t)>2$.  In this case,
\[
\frac{1}{t^{A}\varepsilon_s^{1/2}\varepsilon_t^{1/2}}
\leq
\frac{1}{t^{A}}\left(\sum_{r=a}^{b}\psi(n_r)\right)^{\theta},
\]
and so
\begin{equation}
\begin{aligned} \label{sum_const_leq_linear_case1}
	\mathop{\sum\sum}\limits_{a\leq s<t\leq b} \frac{1}{t^{A}\varepsilon_s^{1/2}\varepsilon_t^{1/2}} &\leq  \mathop{\sum\sum}\limits_{a\leq s<t\leq b} \frac{1}{t^{A}}\left(\sum_{r=a}^{b}\psi(n_r)\right)^{\theta} \\ &\ll \left(\sum_{r=a}^{b}\psi(n_r)\right)^{\theta}\leq \sum_{r=a}^{b}\psi(n_r).
\end{aligned}
\end{equation}

\medskip

\noindent Case B: $\sum\limits_{t=a}^{b}\psi(n_t)\leq 2$. In this case, we have $\varepsilon_t=2^{-\theta}$ ($a\leq t\leq b$), so
\[
\frac{1}{t^{A}\varepsilon_s^{1/2}\varepsilon_t^{1/2}}\ll\frac{1}{t^{A}}.
\]
Applying Corollary~\ref{cor_min_psi} with $\tau = A/2$, we may assume that
\[
\frac{1}{t^{A}}\leq\frac{\psi(n_t)}{t^{A/2}}
\qquad (t \in \bN),
\]
whence
\begin{equation} 
\begin{aligned}
	\label{sum_const_leq_linear_case2}
	\mathop{\sum\sum}\limits_{a\leq s<t\leq b} \frac{1}{t^{A}\varepsilon_s^{1/2}
		\varepsilon_t^{1/2}} &\ll \mathop{\sum\sum}\limits_{a\leq s<t\leq b} {\psi(n_t)}{t^{-A/2}} \\ &\leq \mathop{\sum\sum}\limits_{a\leq s<t\leq b} {\psi(n_t)}{s^{-A/2}}  \ll \sum_{t=a}^b\psi(n_t).
\end{aligned}
\end{equation}
It follows from~\eqref{sum_const_leq_linear_case1} and~\eqref{sum_const_leq_linear_case2} that, in any case,
\begin{equation} \label{sum_const_leq_linear}
\mathop{\sum\sum}\limits_{a\leq s<t\leq b} \frac{1}{t^{A}\varepsilon_s^{1/2}\varepsilon_t^{1/2}}\ll \sum_{t=a}^b\psi(n_t).
\end{equation}

We now proceed to estimate the first term on the right-hand side of~\eqref{ub_sum_mu}. To this end, we will make use of~\eqref{main_q_Harman_ub2}.  With this in mind, first note that by definition $t \mapsto \eps_t$ is non-increasing, so 
\[
\eps_t \psi(n_s) \psi(n_t)
\leq
\eps_s \psi(n_s) \psi(n_t)
\]
whenever $s < t$. Moreover, we have $\eps_t \le 1$ for all $t\in\bN$, and so 
\[
\varepsilon_s
\varepsilon_t
\psi(n_s) \psi(n_t) \le \varepsilon_s
\psi(n_s)\psi(n_t).
\]
Thus
\begin{equation} \label{sum_eps_psipsi_developed}
4(1+\varepsilon_s)(1+\varepsilon_t)\psi(n_s)\psi(n_t) \, \leq \, 4\psi(n_s)\psi(n_t) + 12 \varepsilon_s \psi(n_s) \psi(n_t).
\end{equation}

We estimate the sum of the second term on the right-hand side of \eqref{sum_eps_psipsi_developed} by considering two cases, A and B as above.

\noindent Case A: $\sum\limits_{t=a}^{b}\psi(n_t)\geq 2 $. In this case, we have
\begin{align}
\notag
&\mathop{\sum\sum}\limits_{a\leq s<t\leq b} \varepsilon_s\psi(n_s)\psi(n_t) \\
\notag
&= \mathop{\sum\sum}\limits_{a\leq s<t\leq b} \psi(n_s)\psi(n_t)\ \min\Bigg(2^{-\theta}, \left(  \sum_{r=a}^{s}\psi(n_r)
\right)^{-\theta} \Bigg) \\ 
\notag
&\leq \
\max\Bigg(2^{1-\theta},\left(\sum_{t=a}^b\psi(n_t) \right)^{1-\theta}\Bigg)\mathop{\sum\sum}\limits_{a\leq s<t\leq b}\frac{\psi(n_s)\psi(n_t)}{\max\left(2, \sum_{r=a}^{s}\psi(n_r) \right)} \\
\label{ub_sum_eps_psipsi} &\leq 
\left(\sum_{t=a}^b\psi(n_t) \right)^{1-\theta}\sum_{a\leq t\leq b}\psi(n_t)\sum_{a\leq s\leq b} \frac{\psi(n_s)}{\max\left(2, \sum_{r=a}^{s}\psi(n_r) \right)}.
\end{align}
Applying \cite[Lemma~D2]{PVZZ2022} with 
\[
\gamma = 2,
\qquad
s_k = \psi(n_{k+a-1}),
\]
and $a,b$ replaced by $1,b-a+1$ respectively, gives
\begin{align*}
&\sum_{a\leq s\leq b} \frac{\psi(n_s)}{\max\left(2, \sum_{r=a}^{s}\psi(n_r) \right)} \\
&\qquad <
\frac32 + \frac{1}{2\log\frac32} \left(
\log \left(\sum_{t=a}^b\psi(n_t) \right) - \log \psi(n_a) \right).
\end{align*}
However, since we must be in Case (ii) or (iii) within the proof of that lemma, and since $\sum_{t=a}^b \psi(n_t) \ge 2$, we obtain
\[
\sum_{a\leq s\leq b} \frac{\psi(n_s)}{\max\left(2, \sum_{r=a}^{s}\psi(n_r) \right)} <
\frac32 + \frac{1}{2\log\frac32}\log \left(\sum_{t=a}^b\psi(n_t) \right).
\]
This, together with~\eqref{ub_sum_eps_psipsi}, implies that
\begin{align} 
\notag
&\mathop{\sum\sum}\limits_{a\leq s<t\leq b} \varepsilon_s\psi(n_s)\psi(n_t) \\
\notag
&\leq \left(\sum_{t=a}^b\psi(n_t) \right)^{1-\theta}\left(\sum_{t=a}^b\psi(n_t) \right)\left(\frac32 + \frac{1}{2\log\frac32} \log \left(\sum_{t=a}^b\psi(n_t) \right)\right)\\[2ex]
\label{ub_sum_eps_psipsi_2}
&\ll \ \left( \sum_{t=a}^b\psi(n_t)
\right)^{2-\theta}
\log^+ \left( \sum_{t=a}^b\psi(n_t) \right).
\end{align}

\medskip

\noindent Case B:  $\sum\limits_{t=a}^{b}\psi(n_t)<2$. In this case, we have
\begin{align}
\notag
\mathop{\sum\sum}
\limits_{a\leq s<t\leq b} \varepsilon_s\psi(n_s)\psi(n_t)
&\leq \mathop{\sum\sum}
\limits_{a \leq s < t \leq b} \psi(n_s) \psi(n_t)
\le \left(\sum_{t=a}^b\psi(n_t)\right)^2 \\
\label{ub_sum_eps_psipsi_simple}
&< 2\sum_{t=a}^b\psi(n_t).
\end{align}

Combining  \eqref{sum_eps_psipsi_developed},
\eqref{ub_sum_eps_psipsi_2} and \eqref{ub_sum_eps_psipsi_simple} gives
\begin{equation}
\begin{aligned} \label{last_inequality}
	&\mathop{\sum\sum}\limits_{a\leq s<t\leq b}  (1+\varepsilon_s)(1+\varepsilon_t)\psi(n_s)\psi(n_t) \\ 
	&\leq  \left( \sum_{t=a}^{b}\psi(n_t)\right)^2 
	+  O\left(\left(\sum_{t=a}^b\psi(n_t) \right)^{2-\theta}\log^+ \left(\sum_{t=a}^b\psi(n_t) \right)\right)
	\\ &\qquad + O\left( 
	\sum_{t=a}^b \psi(n_t)
	\right).
\end{aligned}
\end{equation}
The result follows on combining
\eqref{ub_sum_mu}, \eqref{ub_first_two_terms}, \eqref{sum_const_leq_linear} and~\eqref{last_inequality}.
\end{proof}

\bigskip

\begin{proof} [Proof of Proposition~\ref{indieMainThm}]
By employing Proposition~\ref{indiethm_general} with $\theta=2/3$, we find that it is enough to provide upper bounds for the sums
\begin{equation*}
\mathop{\sum\sum}_{a\leq s<t\leq b} (n_s,n_t)\min\Big(\frac{\psi(n_s)}{n_s},\frac{\psi(n_t)}{n_t}\Big)
\end{equation*}
and
\begin{equation} \label{sum_Tst}
\mathop{\sum\sum}\limits_{a\leq s<t\leq b} T(s,t),
\end{equation}
where $T(s,t)$ is given by \eqref{def_S_prime}.
We now carry this out.

First, using the fact that $(n_t)_{t\in\bN}$ is a lacunary sequence, we find
\begin{equation*} \label{ub_gcd}
\begin{aligned}
	\mathop{\sum\sum}_{a\leq s<t\leq b} (n_s,n_t)\min\Big(\frac{\psi(n_s)}{n_s},\frac{\psi(n_t)}{n_t}\Big)
	&\leq \mathop{\sum\sum}_{a\leq s<t\leq b} n_s \frac{\psi(n_t)}{n_t}\\
	=\sum_{t=a+1}^b \psi(n_t)\sum_{a\leq s < t}\frac{n_s}{n_t}
	&\ll \sum_{t=a+1}^b \psi(n_t).
\end{aligned}
\end{equation*}
Next, we are going to find a suitable upper bound for the sum~\eqref{sum_Tst}.

Recall that $T(s,t)$ is the sum of terms 
\[
\what{W}_{n_s,\varepsilon_s}^{+}(kn_s,-k)\what{W}_{n_t,\varepsilon_t}^{+}(\ell n_t,-\ell)  \what{\nu}\left( -kn_s-\ell n_t,k+\ell\right),
\]
over non-zero integers $k,\ell$ for which
\begin{equation} \label{condition_T}
|k n_s-\ell n_t|< n_s^{1/2},
\end{equation}
where $\eps_a, \ldots, \eps_b$ are given by~\eqref{def_varepsilon_t}. By dividing both sides of~\eqref{condition_T} by $n_s$, we get
\begin{equation} \label{index_restriction_first_boundsv}
|k -\ell n_t/n_s|< 1.
\end{equation}

Hence, if $k$ and $\ell$ satisfy \eqref{index_restriction_first_boundsv} then they must have the same sign. Further, for each fixed integer $\ell$, there exists a set $S_{\ell}$ of at most two non-zero integers  $k$  satisfying \eqref{index_restriction_first_boundsv}.
On using the trivial bound  $|\what{\nu}(k_1,k_2)| \leq 1 $, it follows that
$$
|T(s,t)| \leq 2 \mathop{\mathop{\sum}_{ k, \ell \in \bN} }_{k \in S_{\ell}} |\what{W}_{n_s,\varepsilon_s}^{+}(k n_s,-k) \what{W}_{n_t,\varepsilon_t}^{+}(\ell n_t,-\ell) |.
$$

Note from \eqref{index_restriction_first_boundsv} and the inequality $\left|\ell n_t/n_s\right| > 1$ that
\begin{equation}
\label{ratio}
\frac12 \frac{n_s}{n_t}  k  \leq  \frac{n_s}{n_t} \max \{1,  k-1 \}  \, \leq \, \ell\, \leq \frac{n_s}{n_t} (k+1)\leq 2\frac{n_s}{n_t} k.
\end{equation}
Using \eqref{W_ub_psi} to bound $\what{W}_{n_t,\eps_t}^{+}(\ell n_t,-\ell)$, and using \eqref{W_ub_psi}, \eqref{W_ub_1s2} and \eqref{ratio}
to bound  
$\what{W}_{n_s,\eps_s}^{+}(k n_s,-k) $, yields
\begin{eqnarray*}
|T(s,t)| &\ll  & \sum_{\ell  \in \N }
\ \min\left(\frac{n_s^2}{n_t^2}\cdot\frac{1}{\ell^2\psi(n_s)\varepsilon_s},\psi(n_s)\right)\psi(n_t)\\[2ex]
&\ll&
\sum_{\ell \leq \frac{n_s}{n_t\psi(n_s)\varepsilon_s^{1/2}}} \hspace{-5mm}
\psi(n_s)\psi(n_t) +
\sum_{\ell > \frac{n_s}{n_t\psi(n_s)\varepsilon_s^{1/2}}} \frac{n_s^2}{n_t^2}\cdot
\frac{\psi(n_t)}{\ell^2\psi(n_s)
	\varepsilon_s}\\[2ex]
&\ll&
\frac{n_s}{n_t}\frac{\psi(n_t)}{\varepsilon_s^{1/2}}  \leq \ \frac{n_s}{n_t}\frac{\psi(n_t)}{\varepsilon_t^{1/2}}. 
\end{eqnarray*}
The last inequality makes use of the fact that 
$t \mapsto \eps_t$ is non-increasing.

By lacunarity, we have 
\[
\sum\limits_{1\leq s<t}n_s/n_t \ll 1.
\]
Consequently, if 
$a,b \in \N$ with 
$a <  b$, then
\[
\mathop{\sum\sum}_{a\leq s < t\leq b} T(s,t)
\ll \sum_{t=a}^b\sum_{s=a}^{t-1}\frac{n_s}{n_t}\frac{\psi(n_t)}{\varepsilon_t^{1/2}} 
\ll \sum_{t=a}^{b}\frac{\psi(n_t)}{\varepsilon_t^{1/2}}\, .
\]
Finally, recalling that the $\varepsilon_t$ are given by
\eqref{def_varepsilon_t}, we obtain the upper bound
\[
\mathop{\sum\sum}_{a\leq s < t\leq b} T(s,t)  \ll
\left(\sum_{t=a}^{b}\psi(n_t)\right)^{1+\theta/2}.
\]
The choice $\theta=2/3$ delivers \eqref{delta_set_frac23}.
\end{proof}

\section{Multiplicative approximation}
\label{MultSection}

Having established Theorem~\ref{GeneralCountingThm}, we now infer Theorem~\ref{GeneralMultThmQ}. %~\ref{MultThm} and~\ref{MultThmQ}.
Given $\alpha \in \mathcal{K}$ and $\gamma\in\mathbb{R}$, we set 
\[ 
C = \sup \left\{ \frac{1}{t} \log q_t(\alpha): t \in \bN \right\} 
> 0.  
\]

Let $(n_t)_{t\in \bN} \subseteq \bN$ be a lacunary sequence such that 
\begin{equation}\label{nt_rate} \|n_t\alpha - \gamma\| \ls \frac{8}{n_t} \qquad (t\gs 1)  \end{equation}
and
\begin{equation} \label{n_t_is_exponential}
8^t<n_t\leq 4e^{6Ct}.
\end{equation}
The existence of such a sequence is shown in \cite[Lemma 2.1]{CZ2021}. We now define an approximating function $\psi:\mathbb{N}\to [0,1]$ by setting $\psi(1)=0$ and 
\[ \psi(n) = \frac{1}{8 \log n} \qquad (n \gs 2). \]
Then, for $N\in\bN$,
\[
\Psi(N)=\sum_{t=1}^N\psi(n_t)\asymp \log N.
\]

By Theorem~\ref{GeneralCountingThm}, there exists $\cB_1 \subseteq \cB$ of Hausdorff dimension $2$ such that, for all $(\beta,\delta)\in \cB_1$ and all $t\in\bN$, we have
\begin{equation} \label{GeneralMultThmQ_asymp_cN}
\cN(\beta,\delta,\psi,t)\asymp\Psi(t)\asymp\log t.
\end{equation}
Given $T\in\bN$, define $t\in\bN$ by
\[
n_t\leq T < n_{t+1}.
\]
By~\eqref{n_t_is_exponential}, we have $T\asymp n_t$ and $t\asymp\log T$. At the same time, by~\eqref{GeneralMultThmQ_asymp_cN}, there exist $\asymp \log t$ indices $1\leq s\leq t$ such that
\begin{equation} \label{secondary_rate_n_s}
\|n_s\beta - \delta\| \ls\frac{1}{8\log n_s}.
\end{equation}
Multiplying~\eqref{nt_rate} and~\eqref{secondary_rate_n_s} we get that, for all such indices $s$, the inequality~\eqref{strong_Littlewood_ie} holds with $n=n_s\leq n_t\leq T$. 
Therefore,
\[
\widetilde{\cN}(\alpha,\beta,\gamma,\delta,T)\gg \log t\asymp \log\log T,
\]
which proves the claim of Theorem~\ref{GeneralMultThmQ}.

\providecommand{\bysame}{\leavevmode\hbox to3em{\hrulefill}\thinspace}


\begin{thebibliography}{50}

\bibitem{BHV}
V. Beresnevich, A. Haynes and S. Velani, \emph{Sums of reciprocals of fractional parts and multiplicative Diophantine approximation},
Mem. Amer. Math. Soc. \textbf{263} (2020), No. 1276.

\bibitem{BV2023}
V. Beresnevich  and S. Velani, \emph{The divergence Borel--Cantelli lemma revisited}, J. Math. Anal. Appl. \textbf{519} (2023), Paper No. 126750, 21 pp.

\bibitem{CT}
S. Chow and N. Technau,
\emph{Dispersion and Littlewood's conjecture}, arXiv:2307.14871.

\bibitem{CZ2021}
S. Chow and A. Zafeiropoulos, \emph{Fully-inhomogeneous multiplicative diophantine approximation of badly approximable numbers}, Mathematika \textbf{67} (2021), 639--646.

\bibitem{Fal2014}
K. Falconer, \emph{Fractal Geometry: Mathematical Foundations and Applications}, John Wiley \& Sons, 2014.

\bibitem{FFK2021}
K. J. Falconer, J. M. Fraser and A. K\"aenm\"aki, \emph{Minkowski dimension for measures}, Proc. Amer. Math. Soc. \textbf{151} (2023), 779--794.

\bibitem{Good}
I. J. Good, \emph{The fractional dimension theory of continued fractions}, Proc. Camb. Phil. Soc. \textbf{37} (1941), 199-228.

\bibitem{harman}
G. Harman, Metric number theory, London Math. Soc. Lecture Note Ser. (N.S.) \textbf{18,} Clarendon Press, Oxford, 1998.

\bibitem{Jar1928}
I. Jarn\'ik, \emph{Zur metrischen Theorie der diophantischen Approximationen}, Proc. Mat. Fyz.
\textbf{36} (1928), 91--106.

\bibitem{JS2016}
T. Jordan and T. Sahlsten, 
\emph{Fourier transforms of Gibbs measures for the Gauss map}, 
Math. Ann. \textbf{364} (2016), 
983--1023.

\bibitem{Kau1980}
R. Kaufman, \emph{Continued fractions and Fourier transforms}, Mathematika \textbf{27} (1980), 
262--267.

\bibitem{KPV2022}
V. Kleptsyn, M. Pollicott and P. Vytnova, \emph{Uniform lower bounds on the dimension of Bernoulli
convolutions}, Adv. Math.
\textbf{395} (2022), 108090.

\bibitem{Mon1994}
H. L. Montgomery, \emph{Ten lectures on the interface between analytic number theory and harmonic analysis},
CBMS Regional Conference Series in Mathematics \textbf{84}, American Mathematical Society, Providence, RI, 1994.

\bibitem{NZM1991}
I. Niven, H. S. Zuckerman and H. L. Montgomery, \emph{An Introduction to the Theory of Numbers}, Fifth edition, Wiley, 1991.

\bibitem{PV2000}
A. D. Pollington and S. L. Velani, 
{\em On a problem in simultaneous diophantine approximation: Littlewood's conjecture},
Acta Math. \textbf{185} (2000),
287--306.

\bibitem{PVZZ2022} 
A. D. Pollington, S. Velani, A. Zafeiropoulos and E. Zorin, \emph{Inhomogeneous diophantine approximation on $M_0$ sets with restricted denominators}, Int. Math. Res. Not. \textbf{2021.}

\bibitem{QR2003}
M. Queff\'elec and O. Ramar\'e, \emph{Analyse de Fourier des fractions continues \`a quotients restreints}, Enseign. Math. (2) \textbf{49} (2003), 335--356.

\bibitem{RS1992}
A. Rockett and P. Sz\"usz, \emph{Continued Fractions}. World Scientific, Singapore, 1992.

\bibitem{Ste1993}
E. M. Stein, \emph{Harmonic analysis: Real-variable Methods, Orthogonality and Oscillatory Integrals}, Princeton University Press, 1993.

\end{thebibliography}
\end{document}